\documentclass[final,onefignum,onetabnum]{siamart220329}

\usepackage{braket,amsfonts}

\usepackage{array}

\usepackage[caption=false]{subfig}

\usepackage{pgfplots}
\pgfplotsset{compat=1.18}

\newsiamthm{claim}{Claim}
\newsiamremark{remark}{Remark}
\newsiamremark{hypothesis}{Hypothesis}
\crefname{hypothesis}{Hypothesis}{Hypotheses}

\usepackage{algorithmic}
\usepackage{algorithm}
\usepackage{comment}

\usepackage{graphicx,epstopdf}
\usepackage{amsopn}

\usepackage{graphics}
\usepackage{amssymb}
\usepackage{amsmath}
\usepackage{xcolor}
\usepackage{todonotes}
\usepackage{mathrsfs}
\usepackage{multicol}
\usepackage{multirow}
\usepackage{float}
\usepackage{booktabs}
\usepackage{bm}
\usepackage{dsfont}
\usepackage{extarrows}
\usepackage{hyperref}
\hypersetup{
  breaklinks,
  colorlinks,
  citecolor={green!50!black},
  urlcolor={red!50!black},
  linkcolor={green!50!black}
}
\usepackage[nameinlink]{cleveref}
\crefname{subsection}{section}{sections}
\Crefname{algorithm}{Algorithm}{Algorithms}

\makeatletter
\def\cahierline{\gdef\@cahierline}

\cahierline{\textcolor{gray}{Cahier du GERAD G-2024-30}}
\def\ps@firstpage{
\def\@oddfoot{}%
  \let\@evenfoot\@oddfoot
  \def\@oddhead{\hfil\normalfont\small\@cahierline}%
  \let\@evenhead\@oddhead 
}

\def\ps@myheadings{%
    \def\@oddfoot{\normalfont\scriptsize\ttfamily\hfil\@cahierline}
    \def\@evenfoot{\normalfont\scriptsize\ttfamily\@cahierline}
    \def\@evenhead{\rlap{\thepage}\hfil\upshape\footnotesize\leftmark\hfil\scshape\hyperlink{contents}{[toc]}}
    \def\@oddhead{\footnotesize\scshape\hyperlink{contents}{[toc]}\hfil{\upshape\footnotesize\rightmark}\hfil\llap{\thepage
}}
    \let\@mkboth\@gobbletwo
    \let\sectionmark\@gobble
    \let\subsectionmark\@gobble
    }
\makeatother

\newcommand{\R}{\mathds{R}}
\newcommand{\C}{\mathds{C}}
\newcommand{\lammin}{\lambda_1}

\Crefname{ALC@unique}{Line}{Lines}

\usepackage[numbers,sort&compress]{natbib}
\makeatletter
\renewcommand\NAT@bibsetnum[1]{\settowidth\labelwidth{\@biblabel{#1}}%
   \setlength{\leftmargin}{\bibindent}\addtolength{\leftmargin}{\dimexpr\labelwidth+\labelsep\relax}%
   \setlength{\itemindent}{-\bibindent}%
   \setlength{\listparindent}{\itemindent}
\setlength{\itemsep}{\bibsep}\setlength{\parsep}{\z@}%
   \ifNAT@openbib
     \addtolength{\leftmargin}{\bibindent}%
     \setlength{\itemindent}{-\bibindent}%
     \setlength{\listparindent}{\itemindent}%
     \setlength{\parsep}{0pt}%
   \fi
}
\makeatother

\usepackage{etoolbox}
\makeatletter
\patchcmd{\NAT@test}{\else \NAT@nm}{\else \NAT@hyper@{\NAT@nm}}{}{}
\makeatother
\newcommand{\smallfrac}[2]{\hbox{\small $\frac{#1}{#2}$}}

\crefname{thm}{Theorem}{Theorems}

\newtheorem{lem}{Lemma}[section]
\crefname{lem}{Lemma}{Lemmas}

\crefname{prop}{Proposition}{Propositions}
\newtheorem{defn}{Definition}[section]
\Crefname{defn}{Definition}{Definitions}
\newtheorem{rem}{Remark}[section]
\crefname{rem}{Remark}{Remarks}

\newtheorem{ex}{Example}

\newcommand{\pmat}[1]{\begin{pmatrix}#1\end{pmatrix}} 

\newcommand{\T}{^T\!}
\newcommand{\Red}[1]{{\color{red} #1}}

\newcommand{\Range}{\mathop{\mathrm{Range}}}
\newcommand{\Null}{\mathop{\mathrm{Null}}}

\newcommand{\spal}{SPAL}
\newcommand{\spalbb}{SPALBB}

\title{An Inexact augmented Lagrangian algorithm\\ for unsymmetric saddle-point systems}

\author{Na Huang%
  \thanks{Department of Applied Mathematics, College of Science,
  China Agricultural University,
    Beijing, China.
    E-mail: hna@cau.edu.cn.
    Research partially supported by National Natural Science Foundation of China (No.\,12001531).
  }
  \and Yu-Hong Dai%
  \thanks{LSEC, Academy of Mathematics and Systems Science,
    Chinese Academy of Sciences, Beijing, China.
    E-mail: dyh@lsec.cc.ac.cn.
  }
  \and Dominique Orban%
  \thanks{GERAD and
    Department of Mathematics and Industrial Engineering,
    Polytechnique Montr\'eal, QC, Canada.
    E-mail: dominique.orban@gerad.ca.
    Research partially supported by an NSERC Discovery Grant.
  }
  \and Michael A. Saunders%
  \thanks{Systems Optimization Laboratory,
    Department of Management Science and Engineering,
    Stanford University, Stanford, CA, USA.
    E-mail: saunders@stanford.edu.
    Version of \today.%
  }
}

\begin{document}

\maketitle

\thispagestyle{firstpage}
\pagestyle{myheadings}

\begin{abstract}
Augmented Lagrangian (AL) methods are a well known class of algorithms for solving constrained optimization problems. They have been extended to the solution of saddle-point systems of linear equations. We study an AL (\spal) algorithm for unsymmetric saddle-point systems and derive convergence and semi-convergence properties, even when the system is singular. At each step, our {\spal} requires the exact solution of a linear system of the same size but with an SPD (2,2) block. To improve efficiency, we introduce an inexact {\spal} algorithm. We establish its convergence properties under reasonable assumptions. Specifically, we use a gradient method, known as the Barzilai-Borwein (BB) method, to solve the linear system at each iteration. We call the result the augmented Lagrangian BB ({\spalbb}) algorithm and study its convergence. Numerical experiments on test problems from Navier-Stokes equations and coupled Stokes-Darcy flow show that {\spalbb} is more robust and efficient than BICGSTAB and GMRES. {\spalbb} often requires the least CPU time, especially on large systems.
\end{abstract}



\section{Introduction}
We consider the unsymmetric saddle-point system
\begin{equation}\label{saddle1}
   \pmat{ G & B \\
         -B^T & 0}
   \pmat{x\\y}
 = \pmat{f\\g},
\end{equation}
where $B \in \R^{n \times m}~(n\ge m)$, and
$G\in \R^{n\times n}$ is positive definite on the nullspace of $B^T$ but may be unsymmetric and/or singular.
Thus, $x\T Gx>0$ for all nonzero $x\in\Null(B\T\,)$. The change of sign in the second block-row of \eqref{saddle1} makes the matrix semipositive real and positive semistable if $G$ is positive semidefinite \cite{benzi2004}. Linear systems like~\eqref{saddle1} arise 
from certain discretizations of Navier-Stokes equations \citep{elman2007}, mixed and mixed-hybrid finite element approximation of the liquid crystal director model \citep{Ramage2013} and coupled Stokes-Darcy flow \citep{Cai2009}, and within interior methods for constrained optimization \citep{ghannad-orban-saunders-2021,Wright1997}. System \eqref{saddle1} is nonsingular if and only if $B$ has full column rank \cite{Benzi2005}. When $B$ corresponds to a discretized gradient operator, as for example in Navier-Stokes equations \citep{elman2007,greif2004}, then $B$ has low column rank and \eqref{saddle1} is singular. 

Iterative methods for solving saddle-point systems have been studied for decades, such as stationary iterations \citep{bai2017regularized,Benzi2005,zulehner2002analysis}, nonlinear inexact Uzawa methods \citep{cheng2000,hu2006nonlinear,lu2010}, nullspace methods \citep{pestana2016null,scott2020null,scott2022null}, Krylov subspace methods \citep{di2021constraint,gould2014,montoison2023gpmr,orban2017}, and preconditioning techniques \citep{Benzi2005,benzi2008some,dollar2010,rozloznik2002krylov}. Some stationary iterative methods and their semi-convergence have been studied for singular cases \citep{cao2016,Zhang2010,zheng2009}.

Let $Q\in \R^{m\times m}$ be symmetric and positive definite (SPD). If we premultiply the second block-row of~\eqref{saddle1} by $-BQ^{-1}$ and add the result to the first block equation, we find that~\eqref{saddle1} is equivalent to
\begin{equation}\label{alasystem}
\pmat{ G+BQ^{-1}B\T & B \\
      -B^T & 0}
\pmat{x\\y}
=\pmat{f-BQ^{-1}g\\g}.
\end{equation}
\citet{golub2003solving} and \citet{greif2004} showed that 
methods based on~\eqref{alasystem} may have advantages. Indeed, even if $G$ is singular or ill-conditioned, the $(1,1)$ block in~\eqref{alasystem} can be made nonsingular, positive definite or well-conditioned with suitable selections of $Q$. When $G$ is symmetric, the symmetric form
$$
T(Q):=\pmat{ G+BQ^{-1}B\T & B \\
      B^T & 0}
$$
of~\eqref{alasystem} is typically preferred.
\citet{golub2003solving} mainly consider the specific case $Q=\gamma I$, where $\gamma>0$ is constant and $I$ is the identity matrix. They provide analytical observations on the spectrum of $T(\gamma I)$ and show that there is a range of values of $\gamma$ that will improve the condition number of $T(\gamma I)$, as well as the condition number of its $(1,1)$ block and the associated Schur complement. In particular, $\gamma=\|B\|^2/\|G\|$ may often force the norm of the added term $\frac{1}{\gamma }BB\T$ to be of the same magnitude as the norm of $G$. \citet{greif2004} experimentally observe that this special choice is typically effective. Apart from the form of \eqref{alasystem}, they also show that when $G$ is symmetric positive semidefinite of nullity $1$, an effective approach to maintaining sparsity is to choose the augmented term as $\tau bb^T$, where $b$ is a known vector not orthogonal to the nullspace of $G$, and $\tau>0$ is a constant that approximately minimizes the condition number of $G+\tau bb^T$.

The approach of replacing~\eqref{saddle1} by~\eqref{alasystem} can be regarded as an augmented Lagrangian ({\spal}) method, also called the method of multipliers \citep{Benzi2005,golub2003solving,greif2004}. For an extensive overview of the augmented Lagrangian approach and its applications, we refer to \citep{birgin2014,bertsekas2014}. \citet{awanou2005trivariate} apply the Uzawa method \citep{arrow1958} to~\eqref{alasystem} with $Q=\gamma I$ and propose the following {\spal} (with $k=0, 1, 2, \dots$ and $y_0$ assumed given):
$$
\left\{
  \begin{array}{l}
(G+\frac1\gamma BB\T\,)x_k=f-\frac1\gamma Bg-By_k,\\
y_{k+1} = y_k+\frac1\gamma(B\T x_k+g).
  \end{array}
\right.
$$
By introducing another parameter $\rho$, \citet{awanou2005} further generalize {\spal} as
\begin{equation}\label{ALG2005}
\left\{
  \begin{array}{l}
(G+\frac1\gamma BQ^{-1}B\T\,)x_k=f-\frac1\gamma BQ^{-1}g-By_k,\\
y_{k+1} = y_k+\frac1\rho Q^{-1}(B\T x_k+g),
  \end{array}
\right.
\end{equation}
and give a first convergence analysis for the case of unsymmetric $G$. They say that the proofs in \cite{glowinski1989} using spectral arguments cannot be extended to the nonsymmetric case. Under the assumptions that $x^TGx\ge0$ for all $x$ and $x^TGx=0$ with $B^Tx=0$ implies $x=0$, they verify convergence by proving that $\|y_{k+1}-y_*\|_Q\le\|y_{k}-y_*\|_Q$ and then $x_k$ converges to $x_*$, where $(x_*,y_*)$ is the exact solution of~\eqref{saddle1}. \citet{awanou2005} also say that their numerical experiments for an inexact Uzawa algorithm applied to~\eqref{alasystem} do not illustrate convergence. However, we have not been able to find their implementation of the inexact version and the numerical results.

We focus here on the inexact {\spal}. Based on a simple splitting of the matrix in~\eqref{saddle1}, we propose a stationary iterative method that is theoretically equivalent to~\eqref{ALG2005} when $\gamma=\rho$. Hence, we also call it {\spal}. We derive its convergence and semi-convergence for $B$ of any rank based on spectral arguments (unlike \cite{awanou2005}) and obtain an explicit range of convergence for the parameter in {\spal}. We allow $G$ here to be indefinite. 
Our {\spal} requires an exact solution of a linear system at each step. To improve efficiency, we propose an inexact {\spal} in which the linear system is solved inexactly. We show that it converges to the solution of~\eqref{saddle1} under reasonable conditions. Gradient methods are a class of simple optimization approaches using the negative gradient of the objective function as a search direction. The Barzilai-Borwein (BB) \citep{Barzilai1988two} method is a gradient method for unconstrained optimization and has proved to be efficient for solving large and sparse unconstrained convex quadratic programming, which is equivalent to solving an SPD linear system. 
When $G$ is unsymmetric positive definite (UPD), the linear system \eqref{a11} in {\spal} is UPD as well. We use the BB method to solve this UPD linear system inexactly. We call the resulting method the augmented Lagrangian BB ({\spalbb}) algorithm and establish its convergence under suitable assumptions. Numerical experiments on linear systems from Navier-Stokes equations and coupled Stokes-Darcy flow show that {\spalbb} often solves problems more efficiently than GMRES \citep{Saad1986} and BICGSTAB \citep{Van1992}.

The paper is organized as follows. In \cref{sec:alg}, we introduce the augmented Lagrangian algorithm. Its convergence and semi-convergence are established in \cref{sec:alg-con1} and \cref{sec:alg-con2}. The inexact {\spal} and its convergence analysis are provided in \cref{sec:inalg}. The augmented Lagrangian BB algorithm is presented in \cref{sec:algbb}. Numerical experiments are reported in \cref{sec:num}. Conclusions appear in \cref{sec:con}.

\subsection*{Notation}
For any $H\in \R^{n\times n}$, we write its inverse, transpose, spectral set, nullspace and range space as $H^{-1}$, $H^T$, $\mathrm{sp}(H)$, $\Null(H)$, and $\Range(H)$. For any $x\in \C^{n}$, we write its conjugate transpose as $x^*$. For symmetric $H$, $\lambda_{\min}(H)$ and $\lambda_{\max}(H)$ denote the minimum and maximum eigenvalues. $\|\cdot\|$ denotes the $2$-norm of a vector or matrix. For an $n\times n$ SPD matrix $G$, $\|x\|_G=\sqrt{\langle Gx,x\rangle}=\|G^{\tfrac{1}{2}}x\|$ for all $x\in \R^n$, and $\|H\|_G=\sup\limits_{x\neq 0}\frac{\|Hx\|_G}{\|x\|_G}=\|G^{\tfrac{1}{2}}HG^{-\tfrac{1}{2}}\|$ for all $H\in \R^{n\times n}$. For simplicity, the column vector $(x\T \ y\T\,)\T$ is written $(x,y)$, $a_+:=\max\{0,a\}$, and $1/0:=+\infty$.

\section{Augmented Lagrangian algorithm}\label{sec:alg}

We present {\spal} for solving the unsymmetric saddle-point system~\eqref{saddle1}. 
Let $Q$ be SPD matrix and $\omega>0$. Since
\begin{equation}\label{a32}
  A := \pmat{ G   & B
           \\-B^T & 0}
     = \pmat{ G   & B
           \\-B^T & \omega Q}
     - \pmat{0 & 0
          \\ 0 & \omega Q},
\end{equation}
the saddle-point system~\eqref{saddle1} is equivalent to
$$
   \pmat{ G   & B
       \\-B^T & \omega Q}
   \pmat{x\\y} = \pmat{f\\\omega Qy+g}.
$$
This suggests \Cref{alg1} for solving system~\eqref{saddle1}.

\Cref{l11} shows that it is always possible to choose $Q$ and $\omega$ such that \eqref{a11} is nonsingular, even if $A$ is singular.

If $G$ is symmetric, \eqref{saddle1} is equivalent to the constrained optimization problem
\begin{equation}\label{cop}
 \min_x \ \tfrac{1}{2}x^TGx-f^Tx 
 \mathrm{\quad s.t.\quad} g+B^Tx=0.
\end{equation}
The $k$-th step of the augmented Lagrangian algorithm for \eqref{cop} solves the subproblem
\begin{equation}
   \min_{x}~\tfrac{1}{2}x^TGx-f^Tx+\frac{1}{2\omega}\left\|g+B^Tx+\omega Qy_k\right\|_{Q^{-1}}^2,
\end{equation}

\begin{algorithm}[t]
\caption{The augmented Lagrangian algorithm {\spal} for solving \eqref{saddle1}}
\label{alg1}
\begin{algorithmic}[1]
\STATE{Given $y_0\in \R^{m}$, $\omega>0$, and SPD $Q \in \R^{m \times m}$, set $k=0$.}
\WHILE{a stopping condition is not satisfied}
\STATE{Compute $(x_{k+1}, y_{k+1})$ according to the iteration
\vspace*{-2pt}
\begin{equation}\label{a11}
   \pmat{ G   & B
      \\ -B^T & \omega Q}
   \pmat{x_{k+1} \\ y_{k+1}} = \pmat{f \\ \omega Qy_{k}+g}.
\end{equation}}
\vspace*{-6pt}
\STATE{Increment $k$ by $1$.}
\ENDWHILE
\end{algorithmic}
\end{algorithm}

\noindent
where $y_k$ is an estimate of the Lagrange multiplier. Its optimal solution $x_{k+1}$ satisfies 
\begin{equation}\label{alite}
    (G+\frac1\omega BQ^{-1}B\T\,)x_{k+1}+By_k=f-\frac1\omega BQ^{-1}g.
\end{equation}
The multiplier is updated as 
\begin{equation}\label{almul}
    y_{k+1}
    =\frac{1}{\omega}Q^{-1}(g+B^Tx_{k+1}+\omega Qy_k)
    =y_k+\frac1\omega Q^{-1}(B^Tx_{k+1}+g).
\end{equation}
Note that \eqref{a11} also gives
\eqref{alite}--\eqref{almul}.
Hence, we also call it the augmented Lagrangian algorithm here. Clearly, \Cref{alg1} is theoretically equivalent to~\eqref{ALG2005} if $\gamma=\rho=\omega$. When $G$ is symmetric, the convergence of {\spal} or its variants has been studied in \citep{glowinski1989}. \citet{awanou2005} first gave convergence results for~\eqref{ALG2005} when $G$ is unsymmetric positive semi-definite but positive definite on $\Null(B^T)$, based on analyzing the error $\|y_{k}-y_*\|_Q$, where $(x_*,y_*)$ is the exact solution of~\eqref{saddle1}. Here we give the convergence analysis of {\spal} in a different way, based on the spectral properties of $T$ in \eqref{a41} below. We derive the explicit range of convergence for $\omega$ and do not require $G$ to be positive semi-definite.


We call $A=M-N$ a splitting if $M$ is nonsingular.
Defining $T=M^{-1}N$, we consider the following iteration scheme for solving $A z=\ell$:
\begin{equation}\label{a31}
  z_{k+1}=T z_k+M ^{-1}\ell.
\end{equation}

First, we show that~\eqref{a32} is a splitting of  $A$ in~\eqref{saddle1}. For convenience, we introduce
\begin{align}
   S_Q = G+\dfrac{1}{\omega}BQ^{-1}B^T,
      &\qquad H = \tfrac{1}{2}(G+G^T),    \label{SqH}
\\ M = \pmat{G & B \\ -B^T & \omega Q},
      &\qquad N = \pmat{0 & 0 \\ 0 & \omega Q}. \label{a51}
\end{align}
Note that $S_Q$ is the Schur complement of $\omega Q$ in $M$.

\begin{lem}\label{l10}
Let $G\in \R^{n\times n}$ be unsymmetric but positive definite on $\Null(B\T\,)$, and
\begin{eqnarray}
    \eta = \inf\limits_{x\notin\Null(B^T)}
    \dfrac{x^THx}{x^TB Q^{-1}B^Tx}.        \label{a26}
\end{eqnarray}
For any SPD $Q\in \R^{m\times m}$, if $0<\omega<1/(-\eta)_+$, then $S_Q$ is positive definite.
\end{lem}

\begin{proof}
Since $G$ is positive definite on $\Null(B\T\,)$, so is $H$. Then for any nonzero $x\in\Null(B^T)$, it holds that
$x^T(H+\tfrac{1}{\omega}BQ^{-1}B^T)x=x^THx>0.$
For any $x\notin \Null(B^T)$, as $\eta>-1/\omega$, we have
$$x^T(H+\frac{1}{\omega}BQ^{-1}B^T)x
=x^THx+\frac{1}{\omega}x^TBQ^{-1}B^Tx
\ge (\eta+\frac{1}{\omega})x^TBQ^{-1}B^Tx>0.$$
Hence $S_Q$ is positive definite because, for any nonzero $x\in \R^n$, $ x^T(S_Q+S_Q^T)x=2x^T(H+\tfrac{1} {\omega}BQ^{-1}B^T)x>0.$
\end{proof}


By \Cref{l10} and some algebraic manipulation, we have the following results.

\begin{lem}\label{l11}
Under the same conditions as in \Cref{l10}, $M$ is nonsingular and
\begin{equation}\label{a35}
  M^{-1}=\pmat{S_Q^{-1} & -\dfrac{1}{\omega}S_Q^{-1}BQ^{-1} \\[8pt]
      \dfrac{1}{\omega}Q^{-1}B^TS_Q^{-1} & \dfrac{1}{\omega}Q^{-1}-\dfrac{1}{\omega^2}Q^{-1}B^TS_Q^{-1}BQ^{-1}}.
\end{equation}
\end{lem}

\begin{lem}\label{lemma:iter}
Under the same conditions as in \Cref{l10}, the iteration matrix of \Cref{alg1} is
\begin{equation}\label{a41}
 T=M^{-1}N
 = \smash[t]{\pmat{0 & -S_Q^{-1}B
      \\ 0 & I-\dfrac{1}{\omega}Q^{-1}B^TS_Q^{-1}B}}
\end{equation}
and the eigenvalues of $T$ are $0$ with algebraic multiplicity $n$, $1$ with algebraic multiplicity $m-s$, and the remaining $s$ eigenvalues are $\omega\mu/(1+\omega\mu)$, 
where $s$ is the rank of $B$ and $\mu$ is a generalized eigenvalue of $G$ and $BQ^{-1}B^T$ corresponding to the generalized eigenvector $x\notin\Null(B^T)$.
\end{lem}

\begin{proof}
It follows from \eqref{a51} and \eqref{a35} that
\begin{equation*}
 T=\pmat{ G   & B
      \\ -B^T & \omega Q}^{-1}
   \pmat{0 & 0
      \\ 0 & \omega Q}
 = \pmat{0 & -S_Q^{-1}B
      \\ 0 & I-\dfrac{1}{\omega}Q^{-1}B^TS_Q^{-1}B}.
\end{equation*}
Clearly, $T$ has an eigenvalue $0$ with algebraic multiplicity $n$, and the remaining $m$ eigenvalues are $1-\lambda/\omega$, where $\lambda$ is an eigenvalue of $Q^{-1}B^TS_Q^{-1}B$.

Since $S_Q$ is positive definite and $Q$ is SPD, $Q^{-1}B^TS_Q^{-1}B$ is nonsingular when $B$ has full column rank. Thus, $\lambda=0$ if and only if $B$ is column rank-deficient. In this case, $1$ is an eigenvalue of $T$ with algebraic multiplicity $m-s$. 

If $\lambda\neq0$, note that $Q^{-1}B^TS_Q^{-1}B$ and $S_Q^{-1}BQ^{-1}B^T$ possess the same nonzero eigenvalues, and $\lambda$ is also an eigenvalue of $S_Q^{-1}BQ^{-1}B^T$. Then there exists $x\notin\Null(B^T)$ such that 
$S_Q^{-1}BQ^{-1}B^Tx=\lambda x$. Combining with \eqref{SqH} leads to
\begin{equation}\label{geig1}
    Gx =\dfrac{\omega-\lambda}{\omega\lambda} B Q^{-1}B^Tx.
\end{equation}
Hence there exists a generalized eigenvalue $\mu$ of $G$ and $BQ^{-1}B^T$ corresponding to the generalized eigenvector $x\notin\Null(B^T)$ such that $\mu = \tfrac{\omega-\lambda}{\omega\lambda}$, i.e., $\lambda=\tfrac{\omega}{1+\omega\mu}$.
Therefore, we know that the remaining eigenvalues of
$T$ are
$1-\tfrac{1}{1+\omega\mu}=\tfrac{\omega\mu}{1+\omega\mu}.
$
\end{proof}

We should emphasize that \Cref{l10,l11,lemma:iter} hold even if $B$ has low column rank. From \Cref{l11}, we know that $A=M-N $ is a splitting of $A$. Then the convergence analysis of \Cref{alg1} can be based on the spectral properties of $T=M^{-1}N$. In the following, we discuss the convergence of \Cref{alg1} when $B$ does or does not have full column rank, respectively.

\subsection{Convergence analysis when \texorpdfstring{$B$}{B} has full column rank}\label{sec:alg-con1}

In this case, $A$ is nonsingular and the saddle-point system~\eqref{saddle1} has a unique solution.

\begin{theorem}\label{th1}
Suppose $B\in \R^{n\times m}$ has full column rank and $G\in \R^{n\times n}$ is unsymmetric but positive definite on $\Null(B^T)$. For any SPD $Q\in \R^{m\times m}$, let $\eta$ be defined by~\eqref{a26}. If $0<\omega<1/(-2\eta)_+$, then the sequence $\{x_{k},y_k\}$ produced by \Cref{alg1} converges to the unique solution of saddle-point system~\eqref{saddle1}.
\end{theorem}

\begin{proof}
\Cref{alg1} is convergent if and only if the spectral radius of $T$ is less than $1$ \cite[Theorem 4.1]{saad2003}. Note that $0<\omega<1/(-2\eta)_+\le 1/(-\eta)_+$ and the conditions of \Cref{l10} hold. As $B$ has full column rank, it follows from \Cref{lemma:iter} that $1$ is not an eigenvalue of $T$ and then
\begin{equation}\label{a22}
  \rho(T) = \max_{\mu} \dfrac{\omega|\mu|}{|1+\omega\mu|}
          = \max_{\mu}
  \smash[t]{
  \sqrt{\dfrac{(\omega\mu_1)^2+(\omega\mu_2)^2}
              {(1+\omega\mu_1)^2+(\omega\mu_2)^2}}} ,
\end{equation}
where $\mu=\mu_1+{\rm i}\mu_2$ is the generalized eigenvalue of $G$ and $B Q^{-1}B^T$ corresponding to the generalized eigenvector $x\notin\Null(B^T)$. Since $x\notin\Null(B^T)$ and $Q$ is SPD, we have $x^*B Q^{-1}B^Tx>0$. Combining with \eqref{geig1} gives $\mu=\frac{x^*Gx}{x^*B Q^{-1}B^Tx}$. Then
\begin{equation}\label{munu}
  \mu_1 = \dfrac{x^*(G+G^T)x}{2x^*B Q^{-1}B^Tx}
        = \dfrac{x^*Hx}{x^*B Q^{-1}B^Tx}\ge\eta.
\end{equation}
Note that $\eta>-1/(2\omega)$ and $\omega>0$, so that $1+\omega\mu_1\ge1+\omega\eta>1/2$. This together with~\eqref{a22} leads to $\rho(T) < 1$. Therefore, \Cref{alg1} is convergent.
\end{proof}

\begin{rem}\label{rem:derome}
From~\eqref{a22} we see that $\rho(T)$ decreases with $\omega$. This means that the convergence rate of \Cref{alg1} will improve as $\omega$ decreases. In particular, if $\omega=0$ (which means no splitting), $\rho(T)=0$. \Cref{alg1} then reduces to the exact method for problem~\eqref{saddle1}. This is consistent with~\eqref{a11}, i.e., \Cref{alg1} performs only one iteration. 
In addition, since $ \rho(T)\rightarrow0$ as $|\mu|\rightarrow0$, $Q$ should be chosen such that the generalized eigenvalues of $G$ and $B Q^{-1}B^T$ are very close to $0$. Therefore, we can choose $Q$ with very small norm.
\end{rem}

\begin{rem}
If $G$ is semidefinite, we see that $\eta\ge0$. Then \Cref{alg1} is convergent for any $\omega>0$.
\end{rem}

\subsection{Convergence analysis when \texorpdfstring{$B$}{B} is rank-deficient}\label{sec:alg-con2}

In this case, $A$ is singular. We assume that system~\eqref{saddle1} is solvable and show that \Cref{alg1} is semi-convergent. To this end, we introduce some preliminaries on the semi-convergence of iteration scheme~\eqref{a31} for a general linear system $Az=\ell$.

\begin{defn}(\citet[Lemma 6.13]{Berman1994})
\ Iteration~\eqref{a31} is semi-convergent if, for any initial guess $z_0$, the iteration sequence $\{z_k\}$ produced by~\eqref{a31} converges to a solution $z$ of $A z=\ell$ such that
$ z=(I-T)^{D}M ^{-1}\ell+[I-(I-T)^D(I-T)]z_0, $
where $(I-T)^D$ denotes the Drazin inverse \citep{Campbell1979} of $I-T$.
\end{defn}

\begin{lem}[{\protect \citealp[Theorem 6.19]{Berman1994}}]
\label{l1}
Iteration~\eqref{a31} is semi-convergent if and only if ${\rm index}(I-T)=1$ and $v(T)<1$, where ${\rm index}(I-T)$ is the smallest nonnegative integer $k$ such that the ranks of $(I-T)^k$ and $(I-T)^{k+1}$ are equal, and
$ v(T) = \max\{ |\lambda|:
                 ~\lambda\in{\rm sp}(T),~\lambda\ne1 \}
$
is called the pseudo-spectral radius of $T$.
\end{lem}

\begin{lem}[{\protect \citealp[Theorem 2.5]{Zhang2010}}]
\label{l2}
${\rm index}(I-T)=1$ holds if and only if, for all $0\neq w\in \Range(A)$, $w\notin \Null({A M}^{-1})$, i.e., $\Range(A)\cap\Null({A M}^{-1})=\{0\}$.
\end{lem}

In the following, we analyze the semi-convergence property for \Cref{alg1}. By \Cref{l1}, first, we need to show ${\rm index}(I-T)=1$.

\begin{theorem}\label{th2}
Suppose $B\in \R^{n\times m}$ is rank-deficient and $G\in \R^{n\times n}$ is unsymmetric but positive definite on $\Null(B^T)$. For any SPD $Q\in \R^{m\times m}$, let $\eta$ be defined by~\eqref{a26}. If $0<\omega<1/(-\eta)_+$, then ${\rm index}(I-T)=1$.
\end{theorem}

\begin{proof}
Suppose $0 \neq w\in \Range(A)$. Then there is $v=(v_1,v_2)\in \R^{n+m}$ such that
\begin{equation}\label{a36}
  w = Av = \pmat{G & B \\
                 -B^T & 0} \pmat{v_1 \\ v_2}
    = \pmat{G v_1+Bv_2 \\ -B^Tv_1}
  \ne 0.
\end{equation}
By~\eqref{a35}, we have
\begin{eqnarray}
   {A M}^{-1}w &=&
   \smash[t]{
   \pmat{I & 0
      \\ -B^TS_Q^{-1} & \dfrac{1}{\omega}B^T S_Q^{-1}BQ^{-1}}
   \pmat{G v_1+Bv_2 \\ -B^Tv_1}}   \nonumber
\\ &=&\pmat{G v_1+Bv_2
    \\ -B^TS_Q^{-1}(G v_1+Bv_2) -
       \dfrac{1}{\omega} B^T S_Q^{-1}BQ^{-1}B^Tv_1}.  \label{a34}
\end{eqnarray}

If $G v_1+Bv_2\neq0$, clearly, ${A M}^{-1}w\neq 0$, which shows that $w\notin \Null({A M}^{-1})$.

If $G v_1+Bv_2=0$, it follows from~\eqref{a36} that $B^Tv_1\neq0$ and \eqref{a34} yields
\begin{equation}\label{a37}
   {A M}^{-1}w =
   \pmat{ 0 \\ -\dfrac{1}{\omega}B^T S_Q^{-1}BQ^{-1}B^Tv_1}.
\end{equation}
Note that $Q$ is SPD and $B^Tv_1\neq0$, so that $BQ^{-1}B^Tv_1\neq 0$. 
Then we would have
$$B^T S_Q^{-1}BQ^{-1}B^Tv_1\neq0.$$
Indeed, if $B^T S_Q^{-1}BQ^{-1}B^Tv_1=0$, clearly $v_1^TBQ^{-1}B^T S_Q^{-1}BQ^{-1}B^Tv_1=0$. Since $S_Q$ is positive definite
, $S_Q^{-1}$ is also positive definite, which leads to $BQ^{-1}B^Tv_1=0$. This is a contradiction. Therefore, we still get $w\notin \Null({A M}^{-1})$ by~\eqref{a37}.
Summing up, for any $0\neq w\in \Range(A)$, $w\notin \Null({A M}^{-1})$.
The result follows from \Cref{l2}.
\end{proof}

Next, we show that $v(T)<1$.

\begin{theorem}\label{th3}
Suppose $B\in \R^{n\times m}$ is rank-deficient and $G\in \R^{n\times n}$ is unsymmetric but positive definite on $\Null(B\T\,)$. For any SPD $Q\in \R^{m\times m}$, let $\eta$ be defined by~\eqref{a26}. If $0<\omega<1/(-2\eta)_+$, then $v(T)<1$.
\end{theorem}

\begin{proof}
Since $0<\omega<1/(-2\eta)_+\le1/(-\eta)_+$, the conditions of \Cref{l10} hold. 
Note the definition of the pseudo-spectral radius in \Cref{l1}. From~\Cref{lemma:iter},
\begin{eqnarray*}
v(T)
=\max_{\mu} \dfrac{\omega|\mu|}{|1+\omega\mu|}
=\max_{\mu}  \sqrt{\dfrac{(\omega\mu_1)^2+(\omega\mu_2)^2 }{(1+\omega\mu_1)^2+(\omega\mu_2)^2}},
\end{eqnarray*}
where $\mu=\mu_1+{\rm i}\mu_2$ is the generalized eigenvalue of $G$ and $B Q^{-1}B^T$ that corresponds to the generalized eigenvector $x\notin\Null(B^T)$. By~\eqref{munu}, $\omega>0$ and $\eta>-1/(2\omega)$, we have $1+2\omega\mu_1\ge1+2\omega\eta>0$, giving $v(T)<1$.
\end{proof}

Combining \Cref{l1} with \Cref{th2,th3} and $1/(-2\eta)_+<1/(-\eta)_+$, we get the following convergence result.

\begin{theorem}\label{th4}
Suppose $B\in \R^{n\times m}$ is rank-deficient, and $G\in \R^{n\times n}$ is unsymmetric but positive definite on $\Null(B\T\,)$. For any SPD $Q\in \R^{m\times m}$, let $\eta$ be defined by~\eqref{a26}. If $0<\omega<1/(-2\eta)_+$, then the sequence $\{x_{k},y_k\}$ produced by \Cref{alg1} is semi-convergent to a solution of the singular saddle-point system~\eqref{saddle1}.
\end{theorem}

\section{Inexact augmented Lagrangian algorithm}\label{sec:inalg}

In this section, we develop and analyze inexact {\spal} to solve~\eqref{saddle1}. 
Let $\ell=(f,g)$, $z_k=(x_k,y_k)$, and $r_k=A z_k-\ell$. It follows from \eqref{a31} and $A=M-N$ that \Cref{alg1} is equivalent to
\begin{equation}\label{iter-r}
    z_{k+1}=M ^{-1}N z_k+M ^{-1}\ell
    =M ^{-1}(M-A) z_k+M ^{-1}\ell
    =z_k-M ^{-1}r_k,
\end{equation}
where $M$ and $N $ are defined in~\eqref{a51}. To describe the inexact version of \Cref{alg1}, as done in \cite{hu2006nonlinear}, we introduce a nonlinear mapping $\Psi: \R^{n+m}\longrightarrow  \R^{n+m}$
such that for any given $r\in  \R^{n+m}$, $\Psi(r)$ approximates the solution $\Delta z$ of
   $M\Delta z=r$ in that
\begin{equation}\label{a52}
  \|r-M \Psi(r)\|_*\leq \delta \|r\|_*
\end{equation}
for some $\delta\in[0,1)$ and some norm $\|\cdot\|_*$. We obtain the inexact augmented Lagrangian algorithm of \Cref{inalg1}, where the main idea is to approximate $M ^{-1}r_k$ in~\eqref{iter-r}.

\begin{algorithm}
\caption{Inexact augmented Lagrangian algorithm}
\label{inalg1}
\begin{algorithmic}[1]
\STATE{Given $z_0=(x_0,y_0)\in \R^{n+m}$, $\omega>0$, $0\le\delta<1$ and SPD $Q$, set $k=0$.}
\WHILE{a stopping condition is not satisfied}
\STATE{Compute $r_k=A z_k-\ell$.}
\STATE{Compute $\Psi(r_k)\approx M ^{-1}r_k$ satisfying \eqref{a52}.}
\STATE{Compute $z_{k+1}=z_k-\Psi(r_k)$.}
\STATE{Increment $k$ by $1$.}
\ENDWHILE
\end{algorithmic}
\end{algorithm}

In our convergence analysis we use $\|\cdot\|_P$ in~\eqref{a52}, where $P_{\beta}=
\smash[b]{\pmat{I & 0 \\ 0 & \beta Q^{-1}}}$ is SPD and $\beta>0$ is an arbitrary constant. By \Cref{inalg1},
\begin{eqnarray}
r_{k+1} &=& A z_{k+1}-\ell
          = A(z_k-\Psi(r_k))-\ell
          = r_k-A \Psi(r_k) \nonumber
\\ &=& (I-A M ^{-1})r_k+A M ^{-1}(r_k-M  \Psi(r_k) ) \nonumber
\\ &=& N M ^{-1}r_k+(I-N M ^{-1})(r_k-M  \Psi(r_k) ).
\label{recursion-r}
\end{eqnarray}
Likewise, we discuss the convergence of \Cref{inalg1} when $B$ does or does not have full column rank, respectively.

\subsection{Convergence analysis when \texorpdfstring{$B$}{B} has full column rank}\label{sec:inalg-con1}

Note that $P_{\beta}$ is SPD, and~\eqref{recursion-r} gives
\begin{eqnarray*}
  P_{\beta}^{\tfrac{1}{2}}r_{k+1} =P_{\beta}^{\tfrac{1}{2}}N M ^{-1}P_{\beta}^{-\tfrac{1}{2}}P_{\beta}^{\tfrac{1}{2}}r_k+P_{\beta}^{\tfrac{1}{2}}(I-N M ^{-1})P_{\beta}^{-\tfrac{1}{2}}P_{\beta}^{\tfrac{1}{2}}(r_k-M  \Psi(r_k) ).
\end{eqnarray*}
This along with~\eqref{a52} yields
\begin{eqnarray}
  \|r_{k+1}\|_{P_{\beta}}&\le&\|P_{\beta}^{\tfrac{1}{2}}N M ^{-1}P_{\beta}^{-\tfrac{1}{2}}\|\|r_k\|_{P_{\beta}}
  +\|P_{\beta}^{\tfrac{1}{2}}(I-N M ^{-1})P_{\beta}^{-\tfrac{1}{2}}\| \| r_k-M  \Psi(r_k) \|_{P_{\beta}}\nonumber\\
  &\le&\big(\|P_{\beta}^{\tfrac{1}{2}}N M ^{-1}P_{\beta}^{-\tfrac{1}{2}}\|
  +\delta\|I-P_{\beta}^{\tfrac{1}{2}}N M ^{-1}P_{\beta}^{-\tfrac{1}{2}}\| \big)  \|r_k\|_{P_{\beta}}\nonumber\\
   &=&\big(\| N M ^{-1} \|_{P_{\beta}}
  +\delta\|I- N M ^{-1} \|_{P_{\beta}} \big)  \|r_k\|_{P_{\beta}}.\label{a53}
\end{eqnarray}
The following result provides sufficient conditions for $\| N M ^{-1} \|_{P_{\beta}}<1$.


\begin{lem}\label{l20}
Suppose $B\in \R^{n\times m}$ has full column rank and $G\in \R^{n\times n}$ is unsymmetric but positive definite on $\Null(B\T\,)$. For any $\beta>0$ and SPD $Q\in \R^{m\times m}$, let $\eta$ be defined by~\eqref{a26} and $\lammin$ be the minimum eigenvalue of $2\omega H+BQ^{-1}B^T$. Then, $\lammin>0$ and if $0<\omega<\min\left\{1/(-2\eta)_+,\,\sqrt{\lammin/\beta}\,\right\}$, we have $\| N M ^{-1} \|_{P_{\beta}}<1$.
\end{lem}

\begin{proof}
It follows from $0<\omega<1/(-2\eta)_+\le1/(-\eta)_+$ that $S_Q$ is positive definite. Combining with~\eqref{a51} and~\eqref{a35} leads to
\begin{eqnarray*}
   P_{\beta}^{\tfrac{1}{2}}NM ^{-1}P_{\beta}^{-\tfrac{1}{2}}
&=& P_{\beta}^{\tfrac{1}{2}}\pmat{0 & 0 \\ B^TS_Q^{-1} & I-\dfrac{1}{\omega}B^TS_Q^{-1}B Q^{-1}}P_{\beta}^{-\tfrac{1}{2}}
\\ &=&  \pmat{0 & 0 \\ \sqrt{\beta}Q^{-\tfrac{1}{2}}B^TS_Q^{-1} & I-E}
   =:\widetilde{T},
\end{eqnarray*}
where $E=\frac{1}{\omega}Q^{-\tfrac{1}{2}}B^TS_Q^{-1}B Q^{-\tfrac{1}{2}}$.
This shows that
\begin{equation}\label{norm-spe}
  \| N M ^{-1} \|_{P_{\beta}}=\|P_{\beta}^{\tfrac{1}{2}}N M ^{-1}P_{\beta}^{-\tfrac{1}{2}}\|=\Big(\rho(\widetilde{T}\widetilde{T}^T)\Big)^{\tfrac{1}{2}}.
\end{equation}
By direct calculation and~\eqref{SqH}, we have
\begin{eqnarray}
   \rho\left(\widetilde{T}\widetilde{T}^T\right)
   &=&\rho\left((I-E)(I-E\T\,)+\beta Q^{-\tfrac{1}{2}}B^TS_Q^{-1}S_Q^{-T}B Q^{-\tfrac{1}{2}}\right)\nonumber
\\ &=&\rho\left( I-\tfrac{1}{\omega}Q^{-\tfrac{1}{2}}B^TS_Q^{-1}\Big(S_Q+S_Q^T-\tfrac{1}{\omega}B Q^{-1}B^T -\omega \beta I \Big)S_Q^{-T}B Q^{-\tfrac{1}{2}}
 \right)\nonumber
\\ &=&\rho\left(I-\tfrac{1}{\omega^2}Q^{-\tfrac{1}{2}}B^TS_Q^{-1}\Big(2\omega H+BQ^{-1}B^T-\omega^2\beta I \Big)S_Q^{-T}B Q^{-\tfrac{1}{2}}\right).      \label{spe1}
\end{eqnarray}
Note that $B$ has full column rank and $\omega>0$, and if $2\omega H+BQ^{-1}B^T-\omega^2\beta I$ is SPD, so is $\tfrac{1}{\omega^2}Q^{-\tfrac{1}{2}}B^TS_Q^{-1}\Big(2\omega H+BQ^{-1}B^T-\omega^2\beta I \Big)S_Q^{-T}B Q^{-\tfrac{1}{2}}$. Then all eigenvalues of $\widetilde{T}\widetilde{T}^T$ are less than $1$, i.e., $\rho(\widetilde{T}\widetilde{T}^T)<1$. 
Therefore, in order to prove $\| NM ^{-1} \|_{P_{\beta}}<1$, we just need to find $\omega$ to guarantee that $2\omega H+BQ^{-1}B^T-\omega^2\beta I$ is positive definite. 
Since $H$ is positive definite on $\Null(B^T)$,~\eqref{a26} and $2\omega \eta>-1$ imply $2\omega H+BQ^{-1}B^T$ is positive definite. Thus, $\lammin>0$. Combining with $\omega<\sqrt{\lammin/\beta}$ 
gives the result.  
\end{proof}

\begin{rem}
    The conditions in \Cref{l20} are reasonable. 
    Indeed, for any given $\omega_0\in(0,\,1/(-2\eta)_+)$, $2H+\frac{1}{\omega_0}BQ^{-1}B^T$ is SPD. Then when $0<\omega\le \omega_0$, we have 
    $$
    \lammin
    \ge\lambda_{\min}\left(2\omega H+\tfrac{\omega}{\omega_0} BQ^{-1}B^T\right)
    = \omega\lambda_{\min}\left(2H+\tfrac{1}{\omega_0} B Q^{-1}B^T\right).
    $$
    Then the conditions in \Cref{l20} can be replaced by     $$0<\omega<\min\left\{\omega_0,\,\tfrac{1}{\beta}\lambda_{\min}\left(2H+\tfrac{1}{\omega_0} B Q^{-1}B^T\right)\right\}.$$
    In particular, when $H$ is positive semidefinite, $\eta\ge0$ and $2H+BQ^{-1}B^T$ is SPD. Then we can pick $\omega_0=1$ above and the last condition can be further simplified as
    $$0<\omega<\min\left\{1,\,\tfrac{1}{\beta}\lambda_{\min}(2H+BQ^{-1}B^T)\right\}.$$
\end{rem}


\begin{theorem}\label{th5}
Suppose $B\in \R^{n\times m}$ has full column rank and $G\in \R^{n\times n}$ is unsymmetric but positive definite on $\Null(B\T\,)$. For any $\beta>0$ and SPD $Q\in \R^{m\times m}$, let $\eta$ and $\delta$ be defined by~\eqref{a26} and~\eqref{a52}, and $\lammin > 0$ be the minimum eigenvalue of $2\omega H+BQ^{-1}B^T$. If $\omega$ and $\delta$ satisfy 
$$
0<\omega<\min\left\{\frac{1}{(-2\eta)_+},\,\sqrt{\frac{\lammin}{\beta}}\right\}
\quad\mbox{and}\quad
0\le \delta \le \tfrac{1}{2}\Big(1-\| N M ^{-1} \|_{P_{\beta}}\Big),
$$
then $\{x_{k},y_k\}$ produced by \Cref{inalg1} converges to the unique solution of~\eqref{saddle1}.
\end{theorem}

\begin{proof}
It follows from \Cref{l20} that $\|N M ^{-1} \|_{P_{\beta}}<1$, so that
$\|I- N M ^{-1} \|_{P_{\beta}}\le 1+\|N M ^{-1} \|_{P_{\beta}}<2.$
The result follows from~\eqref{a53} and
\begin{align*}
\| N M ^{-1} \|_{P_{\beta}}+\delta\|I- N M ^{-1} \|_{P_{\beta}}
&\le \| N M ^{-1} \|_{P_{\beta}}+\frac{1-\| N M ^{-1} \|_{P_{\beta}}}{2}\|I- N M ^{-1} \|_{P_{\beta}}\\
&< \| N M ^{-1} \|_{P_{\beta}}+1-\| N M ^{-1} \|_{P_{\beta}}=1.
\tag*{\qed}
\end{align*}
\end{proof}

\begin{rem}
  From \eqref{a53} we have
  $
       \|r_{k}\|_{P_{\beta}}\le\big(\| N M ^{-1} \|_{P_{\beta}}
  +\delta\|I- N M ^{-1} \|_{P_{\beta}} \big)^{k}  \|r_0\|_{P_{\beta}}.
  $
 Hence, based on the conditions of \Cref{th5}, $r_k$ converges to zero linearly. Let $z_*$ be the solution of~\eqref{saddle1}. Then 
 \begin{align*}
   &\|z_{k}-z_*\|_{P_{\beta}} = \|A^{-1}r_k\|_{P_{\beta}}
   = \|P_{\beta}^{\frac12}A^{-1}P_{\beta}^{-\frac12}P_{\beta}^{\frac12}r_k\|
   \le \|P_{\beta}^{\frac12}A^{-1}P_{\beta}^{-\frac12}\|\|P_{\beta}^{\frac12}r_k\|
    \\
   & = \|A^{-1}\|_{P_{\beta}}\|r_k\|_{P_{\beta}}
   \le \|A^{-1}\|_{P_{\beta}}\big(\| N M ^{-1} \|_{P_{\beta}}
  +\delta\|I- N M ^{-1} \|_{P_{\beta}} \big)^{k}  \|r_0\|_{P_{\beta}}
  \\
  & = \|A^{-1}\|_{P_{\beta}}\big(\| N M ^{-1} \|_{P_{\beta}}
  +\delta\|I- N M ^{-1} \|_{P_{\beta}} \big)^{k}  \|A(z_0-z_*)\|_{P_{\beta}}
  \\
  &\le \|A^{-1}\|_{P_{\beta}}\|A\|_{P_{\beta}}\big(\| N M ^{-1} \|_{P_{\beta}}
  +\delta\|I- N M ^{-1} \|_{P_{\beta}} \big)^{k}  \|z_0-z_*\|_{P_{\beta}}.
  \end{align*}
 This implies that $z_k$ converges linearly to $z_*$ under the conditions of \Cref{th5}.
\end{rem}

\begin{rem}\label{rem:convergence}
    If $\beta=\delta$ in \Cref{th5}, since $\omega>0$ and $\delta\ge0$, we know that $\omega<\sqrt{\lammin/\delta}$ holds if and only if $\delta<\lammin/\omega^2$.
    Then the restricted conditions of $\omega$ and $\delta$ in \Cref{th5} can be replaced by
    \begin{equation*}\label{condition}
        0<\omega<\frac{1}{(-2\eta)_+}
\quad\mbox{and}\quad
0\le \delta < \min\left\{\frac{\lammin}{\omega^2},\, \dfrac{1-\| N M ^{-1} \|_{P_{\delta}}}{2}\right\}.
    \end{equation*}
It follows from \eqref{norm-spe} and \eqref{spe1} that    
      \begin{align}     
    &\|NM^{-1}\|_{P_{\delta}}^2
    = \rho(\widetilde{T}\widetilde{T}^T)
    =\rho\Big(I-\tfrac{1}{\omega^2}Q^{-\tfrac{1}{2}}B^TS_Q^{-1}(2\omega H+BQ^{-1}B^T-\delta \omega^2I)S_Q^{-T}B Q^{-\tfrac{1}{2}}\Big)\nonumber\\
    &=\rho\left(I-\tfrac{1}{\omega^2}Q^{-\tfrac{1}{2}}B^TS_Q^{-1}(2\omega H+BQ^{-1}B^T)S_Q^{-T}B Q^{-\tfrac{1}{2}}+\delta Q^{-\tfrac{1}{2}}B^TS_Q^{-1}S_Q^{-T}B Q^{-\tfrac{1}{2}}\right).
    \label{normnm}
    \end{align}
    Note that $\widetilde{T}\widetilde{T}^T$ is symmetric positive semidefinite and $Q^{-\tfrac{1}{2}}B^TS_Q^{-1}S_Q^{-T}B Q^{-\tfrac{1}{2}}$ is SPD, $\|NM^{-1}\|_{P_{\delta}}$ increases with $\delta$, and 
$$
\lim_{\delta\rightarrow\lammin/\omega^2}\|NM^{-1}\|_{P_{\delta}}=1,
\qquad
\lim_{\delta\rightarrow 0^+}\|NM^{-1}\|_{P_{\delta}}=\sqrt{1-\tilde{\lambda}_1/\omega^2}<1,
$$
where $\tilde{\lambda}_1>0$ is the minimum eigenvalue of $Q^{-\tfrac{1}{2}}B^TS_Q^{-1}(2\omega H+BQ^{-1}B^T)S_Q^{-T}B Q^{-\tfrac{1}{2}}$. Then there exists $\delta>0$ such that $\|NM^{-1}\|_{P_{\delta}}<1$. Therefore, for any given $0<\omega<1/(-2\eta)_+$, \Cref{inalg1} is convergent for sufficiently small $\delta$. Moreover, the larger $\omega$ is, the smaller $\delta$ should be. Therefore, a practical selection of $\delta$ could be a sequence $\{\delta_k\}$ such that $\delta_k\rightarrow0$ as $k\rightarrow\infty$.
\end{rem}

\begin{rem}\label{rem:positivecase}
When $G$ is positive semidefinite, \eqref{a26} yields $\eta\ge0$. It leads to $(-2\eta)_+=0$. In this case, the sufficient conditions in \Cref{th5} can be replaced by $0<\omega<\min\sqrt{\lammin/\beta}$ and $0\le \delta \le \tfrac{1}{2}\left(1-\| N M ^{-1} \|_{P_{\beta}}\right)$.
Furthermore, from \Cref{rem:convergence} we know that the restrictions also can be replaced by $\omega>0$ and $0\le \delta < \min\left\{\tfrac{\lammin}{\omega^2},\, \tfrac{1-\| N M ^{-1} \|_{P_{\delta}}}{2}\right\}$.
This implies that when $G$ is positive semidefinite, for any $\omega>0$, \Cref{inalg1} is convergent for sufficiently small $\delta$.
\end{rem}

\subsection{Convergence analysis when \texorpdfstring{$B$}{B} is rank-deficient}\label{sec:inalg-con2}

Assume that the rank of $B$ is $s$ and $0<s<m$. Let $B=U\pmat{\Sigma&0} V^T$ be the singular value decomposition (SVD), where $n\times n$ $U$ and $m\times m$ $V$ are orthogonal matrices, $\Sigma=\pmat{\Sigma_s\\0}\in\R^{n\times s}$ has full column rank, and $\Sigma_s={\rm diag}\{\sigma_1,\sigma_2,\ldots,\sigma_s\}$ with all $\sigma_j>0$ contains the singular values of $B$. Let $Q_1\in\R^{s\times s}$ and $Q_2\in\R^{(m-s)\times(m-s)}$ be SPD, and
\begin{align}
&Q=V\pmat{ Q_1 & 0
        \\ 0  & Q_2 }V\T,  \qquad~
\widetilde{D}=\pmat{ U & 0
                 \\  0 & V }, \qquad\quad~
\widetilde{P}_{\beta}=\pmat{ I & 0
                 \\  0 & \beta Q_1^{-1} },                
\label{def:qp} 
\\   
&\widetilde{A}=\pmat{ U^TGU  & \Sigma
                 \\ -\Sigma & 0 }, \qquad
\widetilde{M}=\pmat{ U^TGU & \Sigma
                 \\ -\Sigma & \omega Q_1 }, \qquad
\widetilde{N}=\pmat{0 & 0
                \\  0 & \omega Q_1}.\label{def:mn}
\end{align}
Let $\widetilde{r}_k=\widetilde{D}\T r_k=\pmat{\widetilde{r}_k^{a},\,\widetilde{r}_k^{b}}$, $\widetilde{\Psi}(r_k)=\widetilde{D}\T \Psi(r_k)=
\pmat{\widetilde{\Psi}^a(r_k),\,\widetilde{\Psi}^b(r_k)}$ with $\widetilde{r}_k^{a},\,\widetilde{\Psi}^a(r_k)\in\R^{n+s}$. It follows from~\eqref{a32}, \eqref{a51}, \eqref{def:qp} and~\eqref{def:mn} that
\begin{align}
\widetilde{D}\T A\widetilde{D}
&=\pmat{ U\T & 0
    \\  0   & V\T }
 \pmat{ G   & B
     \\-B^T & 0}
 \pmat{U & 0
    \\ 0 & V }
=\pmat{ U\T GU   & U\T BV
     \\-V\T B^TU & 0}\nonumber\\
&=\pmat{ U\T GU   & \Sigma  & 0
    \\ -\Sigma\T &    0    & 0
    \\     0     &    0    & 0 }
=:\pmat{ \widetilde{A}   & 0 \\  0  & 0 },\label{def:tilA}
\\
\widetilde{D}\T M\widetilde{D}
&=\pmat{ U\T & 0
    \\  0   & V\T }
 \pmat{G & B \\ -B^T & \omega Q}
 \pmat{U & 0
    \\ 0 & V }
=\pmat{ U\T GU   & U\T BV
     \\-V\T B^TU & \omega V\T QV}\nonumber\\
&=\pmat{ U\T GU   & \Sigma     & 0
    \\ -\Sigma\T & \omega Q_1 & 0
    \\     0     &    0       & \omega Q_2 }
=:\pmat{ \widetilde{M}   & 0 \\  0  & \omega Q_2 },\label{def:tilM}
\\
\widetilde{D}\T N\widetilde{D}
&=\pmat{ U\T & 0
    \\  0   & V\T }
 \pmat{0 & 0 \\ 0 & \omega Q}
 \pmat{U & 0
    \\ 0 & V }
=\pmat{ 0   & 0
     \\ 0   & \omega V\T Q V}\nonumber\\
&=\pmat{ 0   &     0      & 0
    \\  0   & \omega Q_1 & 0
    \\  0   &    0       & \omega Q_2 }
=:\pmat{ \widetilde{N}   & 0 \\  0  & \omega Q_2 }.\label{def:tilN}
\end{align}
Based on the above notations, we have the following results.

\begin{lem}\label{lem:rkbzero}
   Suppose $B\in\R^{n\times m}$ is rank-deficient with rank $s$. If~\eqref{saddle1} is solvable, then $\widetilde{r}_k^b=0$ for all $k\ge1$.
\end{lem}

\begin{proof}
    Let $z_*$ be a solution of~\eqref{saddle1}, and let $\widetilde{z}_*=\widetilde{D}\T z_*=\pmat{\widetilde{z}_*^a,\,\widetilde{z}_*^b}$, $\widetilde{z}_k=\widetilde{D}\T z_k=\pmat{\widetilde{z}_k^a,\,\widetilde{z}_k^b}$, and $\widetilde{\ell}=\widetilde{D}\T \ell=\pmat{\widetilde{\ell}^a,\,\widetilde{\ell}^b}$, where $\widetilde{z}_*^a,\,\widetilde{z}_k^a,\,\widetilde{\ell}^a\in\R^{n+s}$. It follows from $Az_*=\ell$ and \eqref{def:tilA} that
$$
\smash[t]{
\widetilde{D}\T A\widetilde{D} \widetilde{z}_*
=\pmat{ \widetilde{A}   & 0 \\  0  & 0 }
\pmat{\widetilde{z}_*^a \\ \widetilde{z}_*^b}
= \pmat{\widetilde{A} \widetilde{z}_*^a \\ 0}
=\pmat{\widetilde{\ell}^a \\ \widetilde{\ell}^b},}
$$
which shows that $\widetilde{\ell}^b=0$. Then we have 
\begin{align*}
\widetilde{r}_k
&=\widetilde{D}\T r_k 
=\widetilde{D}\T (Az_k-\ell)
=\widetilde{D}\T A\widetilde{D}\widetilde{D}\T z_k-\widetilde{D}\T\ell
=\smash[t]{\pmat{\widetilde{A}\widetilde{z}_k^a-\widetilde{\ell}^a \\ -\widetilde{\ell}^b}}
=\pmat{\widetilde{r}_k^a \\0}.
\tag*{\qed}
\end{align*}
\end{proof}

\begin{lem}\label{lem:stop}
   Suppose $B\in\R^{n\times m}$ is rank-deficient with rank $s$. For any $\omega,\,\beta>0$ and SPD $Q_1\in \R^{(n+s)\times (n+s)}$ and $Q_2\in \R^{(m-s)\times (m-s)}$, let $Q$ and $\delta$ be defined by~\eqref{def:qp} and~\eqref{a52}. Then 
   $\|\widetilde{r}^a_k-\widetilde{M}\widetilde{\Psi}^a(r_k)\|_{\widetilde{P}_{\beta}}
   \le \delta\|\widetilde{r}^a_k\|_{\widetilde{P}_{\beta}}$.
\end{lem}

\begin{proof}
    For any $x\in\R^{n+m}$ and $\widetilde{x}=\widetilde{D}^Tx=\pmat{\widetilde{x}^a,\,\widetilde{x}^b}$ with $\widetilde{x}^a\in\R^{n+s}$, since $\widetilde{D}$ is an orthogonal matrix, from \eqref{def:qp} and the definition of $P_{\beta}$ in \cref{sec:inalg}, we have 
\begin{align}
\|x\|_{P_{\beta}}^2
& = x^TP_{\beta}x 
= x^T\widetilde{D}\widetilde{D}^TP_{\beta}\widetilde{D}\widetilde{D}^Tx
=\pmat{(\widetilde{x}^a)^T \,(\widetilde{x}^b)^T}
\pmat{\widetilde{P}_{\beta} & 0 \\ 0 & \beta Q_2^{-1}}\pmat{\widetilde{x}^a \\ \widetilde{x}^b}
\nonumber\\
&=\|\widetilde{x}^a\|_{\widetilde{P}_{\beta}}^2+\|\widetilde{x}^b \|_{\beta Q_2^{-1}}^2.\label{norm:ab}
\end{align}
Note that~\eqref{def:tilM} and~\Cref{lem:rkbzero} give
\begin{align*}
\widetilde{D}^T\left(r_k-M\Psi(r_k)\right)
&=\widetilde{r}_k-\widetilde{D}^T M\widetilde{D}\widetilde{\Psi}(r_k)
= \smash{\pmat{\widetilde{r}^a_k-\widetilde{M}\widetilde{\Psi}^a(r_k) 
\\ -\omega Q_2\widetilde{\Psi}^b(r_k)}}.
\end{align*}
This along with \eqref{norm:ab} leads to
\begin{align*}
&\|r_k-M\Psi(r_k)\|_{P_{\beta}}^2
=\|\widetilde{r}^a_k-\widetilde{M}\widetilde{\Psi}^a(r_k)\|_{\widetilde{P}_{\beta}}^2
+\|\omega Q_2\widetilde{\Psi}^b(r_k)\|_{\beta Q_2^{-1}}^2
\\
&=\|\widetilde{r}^a_k-\widetilde{M}\widetilde{\Psi}^a(r_k)\|_{\widetilde{P}_{\beta}}^2
+\omega^2\beta\|\widetilde{\Psi}^b(r_k)\|_{Q_2}^2.
\end{align*}
Using~\eqref{a52}, \eqref{norm:ab} and $\widetilde{r}^b_k=0$ yields
\begin{align*}
\|\widetilde{r}^a_k-\widetilde{M}\widetilde{\Psi}^a(r_k)\|_{\widetilde{P}_{\beta}}
\le 
\|r_k-M\Psi(r_k)\|_{P_{\beta}}
\le 
\delta \|r_k\|_{P_{\beta}}
=\delta\|\widetilde{r}^a_k\|_{\widetilde{P}_{\beta}}.
\tag*{\qed}
\end{align*}
\end{proof}

We are now ready to establish the convergence theorem for \Cref{inalg1} when $B$ is rank-deficient.

\begin{theorem}\label{th-inexact-defi}
Suppose $B\in\R^{n\times m}$ is rank-deficient with rank $s$ and $G\in \R^{n\times n}$ is unsymmetric but positive definite on $\Null(B\T\,)$. For any $\beta>0$ and SPD $Q_1\in \R^{(n+s)\times (n+s)}$ and $Q_2\in \R^{(m-s)\times (m-s)}$, let $Q$, $\eta$ and $\delta$ be defined by~\eqref{def:qp}, \eqref{a26} and~\eqref{a52}, and $\lammin$ be the minimum eigenvalue of $2\omega H+BQ^{-1}B^T$. If $\omega$ and $\delta$ satisfy 
$$
0<\omega<\min\left\{\frac{1}{(-2\eta)_+},\,\sqrt{\frac{\lammin}{\beta}}\right\}
\quad\mbox{and}\quad
0\le \delta \le \tfrac{1}{2}\Big(1-\| \widetilde{N} \widetilde{M} ^{-1} \|_{\widetilde{P}_{\beta}}\Big),
$$
then $\{x_{k},y_k\}$ produced by \Cref{inalg1} converges to a solution of the singular saddle-point system~\eqref{saddle1}.
\end{theorem}

\begin{proof}
By \Cref{lem:rkbzero}, we just need to prove $\lim\limits_{k\rightarrow0}\widetilde{r}_{k}^{a}=0$. Since $\widetilde{D}$ is an orthogonal matrix, it follows from~\eqref{recursion-r}, \eqref{def:qp}, \eqref{def:tilM} and~\eqref{def:tilN} that
\begin{align*}
&\pmat{\widetilde{r}_{k+1}^{a} \\ \widetilde{r}_{k+1}^{b}} = 
\widetilde{r}_{k+1} =  \widetilde{D}\T r_{k+1} 
=\widetilde{D}\T\left[N M ^{-1}r_k+(I-N M ^{-1})(r_k-M  \Psi(r_k) )\right]
\\
&=\widetilde{D}\T N\widetilde{D} (\widetilde{D}\T M\widetilde{D} )^{-1}\widetilde{D}\T r_k+
\left[I-\widetilde{D}\T N\widetilde{D} (\widetilde{D}\T M\widetilde{D} ) ^{-1}\right]\left(\widetilde{D}\T r_k-\widetilde{D}\T M\widetilde{D} \widetilde{D}\T\Psi(r_k) \right)\\
&=\pmat{ \widetilde{N}\widetilde{M}^{-1}   & 0 \\  0  & I }
\pmat{\widetilde{r}_k^{a} \\ \widetilde{r}_k^{b}}+
\left[I-\pmat{ \widetilde{N}\widetilde{M}^{-1}   & 0 \\  0  & I }\right]
\left[\pmat{\widetilde{r}_k^{a} \\ \widetilde{r}_k^{b}}-
\pmat{ \widetilde{M}   & 0 \\  0  & \omega Q_2 }\pmat{\widetilde{\Psi}^a(r_k) \\ \widetilde{\Psi}^b(r_k)}\right]
\\
&=\pmat{\widetilde{N}\widetilde{M}^{-1}\widetilde{r}_k^{a} + (I-\widetilde{N}\widetilde{M}^{-1})(\widetilde{r}_k^{a}-\widetilde{M}\widetilde{\Psi}^a(r_k))\\
\widetilde{r}_k^{b}}.
\end{align*}
Thus,
$\widetilde{r}_{k+1}^{a}
=\widetilde{N}\widetilde{M}^{-1}\widetilde{r}_k^{a} + (I-\widetilde{N}\widetilde{M}^{-1})(\widetilde{r}_k^{a}-\widetilde{M}\widetilde{\Psi}^a(r_k)).$
Using~\eqref{recursion-r}, \eqref{def:tilA}, \eqref{def:tilM}, \eqref{def:tilN} and~\Cref{lem:stop}, we know that $\widetilde{r}_{k}^{a}$ is the $k$-th residual of \Cref{inalg1} applying to the saddle-point problem $\widetilde{A}\widetilde{z}=\widetilde{\ell}$.

Note that $x\in\Null(\Sigma\T\,)$ if and only if $Ux\in\Null(B\T\,)$ and $U\T G U$ is positive definite on $\Null(\Sigma\T\,)$. With \eqref{a26}, \eqref{def:qp}, and the SVD of $B$, we have
\begin{align*}
&\inf\limits_{x\notin\Null(\Sigma\T)}
\dfrac{x^TU^THUx}{x^T\Sigma Q_1^{-1}\Sigma^Tx}
\xlongequal{\hat{x}=Ux}
\inf\limits_{\hat{x}\notin\Null(B\T)}
\dfrac{\hat{x}^T H \hat{x}}{\hat{x}^TU\Sigma Q_1^{-1} \Sigma^TU^T\hat{x}}\\
&=\inf\limits_{\hat{x}\notin\Null(B\T)}
\dfrac{\hat{x}^TH\hat{x}}{\hat{x}^TU\pmat{\Sigma&0}V^TQ^{-1}V\pmat{\Sigma^T\\0}U^T\hat{x}}
=\inf\limits_{\hat{x}\notin\Null(B\T)}
\dfrac{\hat{x}^TH\hat{x}}{\hat{x}^TBQ^{-1}B^T\hat{x}}
=\eta.
\end{align*}
Since $\omega(U^TGU+U^TG^TU)+\Sigma Q_1^{-1}\Sigma^T$ is similar to $2\omega H+U\Sigma Q_1^{-1}\Sigma^TU^T=2\omega H+BQ^{-1}B^T$ and $\Sigma$ has full 
rank,~\Cref{l20} and~\Cref{th5} imply $\| \widetilde{N} \widetilde{M} ^{-1} \|_{\widetilde{P}_{\beta}}<1$ and hence $\widetilde{r}_{k}^{a}$ converges to zero as $k\rightarrow\infty$. Combining with~$\widetilde{r}_{k}^{b}=0$ concludes.
\end{proof}

Similar to \Cref{rem:convergence,rem:positivecase}, when $B$ is rank-deficient, for any given $0<\omega<1/(-2\eta)_+$, \Cref{inalg1} is still convergent for sufficiently small $\delta\ge0$. Furthermore, when $G$ is positive semidefinite, \Cref{inalg1} is convergent for any $\omega>0$ and sufficiently small $\delta\ge0$.

\subsection{Augmented Lagrangian BB algorithm}\label{sec:algbb}

Gradient-type iterative methods for the unconstrained optimization problem $\min\limits_{z \in \R^{\hat{n}}} \hat{f}(z)$
have the form
\begin{equation}\label{eqitr}
z_{k+1}=z_k-\alpha_kg_k,
\end{equation}
where $\hat{f}: \R^{\hat{n}} \rightarrow \R$ is a sufficiently smooth function, $g_k=\nabla \hat{f}(z_k)$ is the gradient, and $\alpha_k>0$ is a stepsize. Methods of this type differ in their stepsize rules. In 1988, \citet{Barzilai1988two} proposed two choices of $\alpha_k$, usually referred to as the BB method: 
\begin{equation}\label{bb}
\alpha_k^{\rm BB1}=
\smash[t]{
\frac{s_{k-1}^Ts_{k-1}}{s_{k-1}^Td_{k-1}} \quad \textrm{and} \quad \alpha_k^{\rm BB2}=\frac{s_{k-1}^Td_{k-1}}{d_{k-1}^Td_{k-1}},}
\end{equation}
where $s_{k-1}=z_k-z_{k-1}$ and $d_{k-1}=g_k-g_{k-1}$.
The rationale behind these choices is related to viewing the gradient-type methods as quasi-Newton methods, where $\alpha_k$ in~\eqref{eqitr} is replaced by $D_k=\alpha_kI$. This matrix serves as an approximate inverse Hessian. Following the quasi-Newton approach, the stepsize is calculated by forcing either $D_k^{-1}$ (BB1 method) or $D_k$ (BB2 method) to satisfy the secant equation in the least squares sense. The corresponding problems are $\min\limits_{D=\alpha I}~\|D^{-1}s_{k-1}-d_{k-1}\|$ and $\min\limits_{D=\alpha I}~\|s_{k-1}-Dd_{k-1}\|$.

When $\hat{f}(z)$ is a convex quadratic, i.e., $\hat{f}(z) = \tfrac{1}{2}z^T\hat{A}z - \hat{\ell}^Tz$ with $\hat{A}$ SPD, this quadratic programming is equivalent to $\hat{A}z=\hat{\ell}$. In this case, $g_k=\hat{A}z_k-\hat{\ell}=r_k$, 
\begin{equation}\label{sd}
    s_{k-1}=-\alpha_{k-1}r_{k-1}
    \quad\mbox{and}\quad
    d_{k-1}=r_k-r_{k-1}=\hat{A}s_{k-1}=-\alpha_{k-1}\hat{A}r_{k-1}.
\end{equation}
Then the two BB stepsizes~\eqref{bb} can be reformulated as 
$$
\alpha_k^{\rm BB1}=\frac{r_{k-1}^Tr_{k-1}}{r_{k-1}^T\hat{A}r_{k-1}} \quad \textrm{and} \quad \alpha_k^{\rm BB2}=\frac{r_{k-1}^T\hat{A}r_{k-1}}{r_{k-1}^T\hat{A}^T\hat{A}r_{k-1}}.
$$

Global convergence of the BB method for minimizing quadratic forms was established by \citet{raydan1993}, and its R-linear convergence rate was established by \citet{dai2002r}. For general strongly convex functions with Lipschitz gradient, the local convergence of the BB method with R-linear rate was rigorously proved by \citet{dai2006cyclic}. Extensive numerical experiments show that the BB method can solve unconstrained optimization problems efficiently and is considerably superior to the steepest descent method \citep{burdakov2019,raydan1997}.
A variety of modifications and extensions of the BB method have been developed for optimization. 

Several researchers used the BB method to solve UPD linear systems. \citet{dai2005analysis} gave an analysis of the BB1 method for two-by-two unsymmetric linear systems. Under mild conditions, they showed that the convergence rate of the BB1 method is $Q$-superlinear if the 
matrix has a double eigenvalue, but only $R$-superlinear if the matrix has two different real eigenvalues. We find that the BB1 method for solving UPD linear systems could be divergent. Indeed, consider 
$$
\hat{A}z := \pmat{1&2\\-2&1}\pmat{x\\y}=\pmat{0\\0}.
$$
Note that $\hat{A}$ has two complex eigenvalues $1\pm 2{\rm i}$. The conditions in \cite{dai2005analysis} do not hold.
It follows from~\eqref{bb} and~\eqref{sd} that $\alpha_k^{\rm BB1}=(s_{k-1}^Ts_{k-1})/(s_{k-1}^T\hat{A}s_{k-1})=1.$ Then, one BB1 iteration gives
$$
z_{k+1}=z_k-r_k
=
\smash[t]{
\pmat{x_k\\y_k}-\pmat{x_k+2y_k\\-2x_k+y_k}
=\pmat{-2y_k\\2x_k}}.
$$
This leads to $\|z_{k+1}\|^2=8\|z_k\|^2$, which means that the sequence $\{z_k\}$ of the BB1 iterations diverges for any initial $z_0\neq0$. 

For quadratic programming with $\hat{A}$ unsymmetric, the minimal gradient method
\citep{kozjakin1982some,krasnosel1952,saad2003} uses the stepsize
$\alpha_k^{\rm MG}=(r_{k}^T\hat{A}r_{k})/(r_{k}^T\hat{A}^T\hat{A}r_{k})$, which gives an optimal residual in each iteration, namely, 
$$\alpha_k^{\rm MG}
=\arg\min_{\alpha>0}\|\hat{A}(z_{k}-\alpha r_k)-b\|
=\arg\min_{\alpha>0}\|r_k-\alpha \hat{A}r_k\|.$$
Therefore, the minimal gradient method is convergent for solving UPD linear systems. Note that the difference between $\alpha_k^{\rm MG}$ and $\alpha_k^{\rm BB2}$ is that one uses $r_k$ and the other uses $r_{k-1}$. The BB2 method can be regarded as the minimal gradient method with delay \citep{friedlander1998}. Gradient methods with delay significantly improve the performance of gradient methods, see \citep{zou2022delayed} and references therein. Hence, we use the BB2 method to derive the new iterates $x_{k+1}$ and $y_{k+1}$ in \Cref{inalg1} when $G$ is positive definite. Then the augmented Lagrangian BB algorithm for solving \eqref{saddle1} is as in~\Cref{kktbb1}.

\begin{algorithm}
\caption{Augmented Lagrangian BB algorithm, {\spalbb}}
\label{kktbb1}
\begin{algorithmic}[1]
\STATE{Given $z_{-1}=(x_{-1},\,y_{-1}),~z_0=(x_0,\,y_0)\in \R^{n+m}$, $\omega>0$, $0\le\delta<1$, and SPD $Q$, compute $r_0=Mz_0-(f,\,\omega Qy_0+g)$ and set $k=0$.}
\WHILE{a stopping condition is not satisfied}
\STATE{Compute $\ell_k=(f,\,\omega Qy_{k}+g)$.}
\WHILE{$\|r_j-M z_j\|_*> \delta \|r_j\|_*$}
\STATE{Compute $s_{j}=z_{j}-z_{j-1}$.}
\STATE{Compute $d_{j}=M s_{j}$.}
\STATE{Compute $r_{j}=M z_{j}-\ell_k$.}
\STATE{Compute $\alpha_{j}=\frac{s_{j}^Td_{j}}{\|d_{j}\|^2}$.}
\STATE{Compute $z_{{j}+1}=z_{j}-\alpha_{j} r_{j}$.}
\ENDWHILE
\STATE{Increment $k$ by $1$.}
\ENDWHILE
\end{algorithmic}
\end{algorithm}

In the following, we establish the convergence of \Cref{kktbb1}. First, under some assumptions, we show that the BB2 method is convergent for solving a general UPD linear system $\hat{A}z=\hat{\ell}$, where the iterative scheme is
$z_{k+1}=z_{k}-\alpha_k^{\rm BB2}r_{k}$ and 
$r_{k}=\hat{A}z_k-\hat{\ell}.$
For convenience, we introduce 
\begin{align}
&\hat{A}_h=\tfrac{1}{2}(\hat{A}+\hat{A}^T),\qquad W=\hat{A}_h^{-1}\hat{A}\T \hat{A}, \nonumber\\
&\theta_j=\max\left\{1-\frac{2u_j}{\lambda_{\min}(W)}+\frac{|\lambda_j|^2}{\lambda_{\min}(W)^2},\,
 1-\frac{2u_j}{\lambda_{\max}(W)}+\frac{|\lambda_j|^2}{\lambda_{\max}(W)^2}\right\}, \label{thetaj}
\end{align}
where $\lambda_j=u_j+{\rm i}v_j~(1\le j\le n)$ are the eigenvalues of $\hat{A}$. When $\hat{A}$ is UPD, we know that $\hat{A}_h$ is SPD and $u_j>0~(1\le j\le n)$. By direct calculation, for all $1\le j\le n$, $\theta_j<1$ holds by $1-\frac{2u_j}{\lambda_{\min}(W)}+\frac{|\lambda_j|^2}{\lambda_{\min}(W)^2}<1$ and
 $1-\frac{2u_j}{\lambda_{\max}(W)}+\frac{|\lambda_j|^2}{\lambda_{\max}(W)^2}<1$, which are equivalent to
\begin{equation}\label{con1}
 \smash[t]{
 \max_{1\le j\le n}\frac{|\lambda_j|^2}{u_j}<2\lambda_{\min}(W).}
\end{equation}
We are now ready to study the convergence of the BB2 method.

\begin{theorem}\label{thbb}
Suppose $\hat{A}\in \R^{\hat{n}\times \hat{n}}$ is UPD. If its $n$ eigenvalues $\lambda_j=u_j+{\rm i}v_j~(1\le j\le n)$ satisfy \eqref{con1}, then the sequence $\{z_{k}\}$ produced by the BB2 method converges to the unique solution of $\hat{A}z=\hat{\ell}$.
\end{theorem}

\begin{proof}
It is well known that the BB method is invariant under unitary transformation of the variables \cite{dai2002r}. By the Schur decomposition, we can assume without loss of generality that $\hat{A}$ is of the form
$$
\pmat{\lambda_1&a_{12}&a_{13}&\cdots&a_{1\hat{n}}\\
0 & \lambda_2&a_{23}&\cdots&a_{2\hat{n}}\\
\vdots&\ddots&\ddots&\ddots&\vdots\\
0&\cdots&0&\lambda_{\hat{n}-1}&a_{\hat{n}-1,\hat{n}}\\
0&\cdots&\cdots&0&\lambda_{\hat{n}}},
$$
where $\lambda_j=u_j+{\rm i}v_j\in\mathbb{C},~j=1,2,\ldots,\hat{n}$. Because $r_{k+1}=\hat{A}z_{k+1}-\hat{\ell}=r_k-\alpha_k^{\rm BB2}\hat{A}r_k$,
\begin{equation}\label{r_k}
  \left\{
    \begin{array}{l}
    r_{k+1}^{(\hat{n})}=r_k^{(\hat{n})}-\alpha_k^{\rm BB2}\lambda_{\hat{n}}r_k^{(\hat{n})},\\[3pt]
    r_{k+1}^{(j)}=r_k^{(j)}-\alpha_k^{\rm BB2}\lambda_jr_k^{(j)} - \alpha_k^{\rm BB2}\sum\limits_{t=j+1}^{\hat{n}}a_{j,t}r_k^{(t)},~j=\hat{n}-1,\ldots,1,
    \end{array}
  \right.
\end{equation}
where $r_k^{(j)}$ is the $j$-th component of $r_k$. 
Note that $\hat{A}_h=\tfrac{1}{2}(\hat{A}+\hat{A}^T)$ and $r_{k-1}^T\hat{A}r_{k-1}=r_{k-1}^T\hat{A}^Tr_{k-1}$, giving
$r_{k-1}^T\hat{A}r_{k-1}=\tfrac{1}{2}\left(r_{k-1}^T\hat{A}r_{k-1}+r_{k-1}^T\hat{A}^Tr_{k-1}\right)=r_{k-1}^T\hat{A}_hr_{k-1}.$
Since $\hat{A}_h$ is SPD, it leads to
$$\alpha_k^{\rm BB2}=\frac{r_{k-1}^T\hat{A}r_{k-1}}{r_{k-1}^T\hat{A}^T\hat{A}r_{k-1}}=\frac{r_{k-1}^T\hat{A}_hr_{k-1}}{r_{k-1}^T\hat{A}^T\hat{A}r_{k-1}}
\xlongequal{\hat{r}=\hat{A}_h^{\frac12}r_{k-1}}
\frac{\hat{r}^T\hat{r}}
{\hat{r}^T\hat{A}_h^{-\frac12}\hat{A}^T\hat{A}\hat{A}_h^{-\frac12}\hat{r}}.$$
By the Courant-Fischer min-max theorem and the fact that $\hat{A}_h^{-\frac12}\hat{A}^T\hat{A}\hat{A}_h^{-\frac12}$ is similar to $W$, we have
\begin{equation}\label{boundBB}
  \smash[t]{
  \frac{1}{\lambda_{\max}(W)}\le\alpha_k^{\rm BB2}\le\frac{1}{\lambda_{\min}(W)}.}
\end{equation}
It follows from $\lambda_j=u_j+{\rm i}v_j$, \eqref{thetaj}, \eqref{boundBB}, and the behavior of the quadratic function for $\alpha_k^{\rm BB2}$ that, for any $j=1,\ldots,\hat{n}$, 
\begin{align}
&\left|1-\alpha_k^{\rm BB2}\lambda_j\right|^2
=\left(1-\alpha_k^{\rm BB2}u_j\right)^2+\left(\alpha_k^{\rm BB2}v_j\right)^2
=1-2\alpha_k^{\rm BB2}u_j+\left(\alpha_k^{\rm BB2}\right)^2\left|\lambda_j\right|^2\nonumber\\
&\le\max\left\{1-\tfrac{2u_j}{\lambda_{\min}(W)}+\tfrac{|\lambda_j|^2}{\lambda_{\min}(W)^2},\,
 1-\tfrac{2u_j}{\lambda_{\max}(W)}+\tfrac{|\lambda_j|^2}{\lambda_{\max}(W)^2}\right\} 
=\theta_j.\label{absal}
\end{align}
Combining with \eqref{con1} and \eqref{r_k} gives
$$
\left|r_{k+1}^{(\hat{n})}\right|
=\left|1-\alpha_k^{\rm BB2}\lambda_{\hat{n}}\right|\,\left|r_{k}^{(\hat{n})}\right|
\le \sqrt{\theta_{\hat{n}}}\left|r_{k}^{(\hat{n})}\right|<\left|r_{k}^{(\hat{n})}\right|.
$$
This implies that $r_{k}^{(\hat{n})}\rightarrow0$ as $k\rightarrow\infty$. For $j=\hat{n}-1,\ldots,1$, by \eqref{r_k} and \eqref{absal}, 
$\left|r_{k+1}^{(j)}\right|
\le
\left|1-\alpha_k^{\rm BB2}\lambda_j\right|\,\left|r_{k}^{(j)}\right| + \alpha_k^{\rm BB2}\left|\sum\limits_{t=j+1}^{\hat{n}}a_{j,t}r_k^{(t)}\right|
\le \sqrt{\theta_j}\left|r_{k}^{(j)}\right| + \alpha_k^{\rm BB2}\left|\sum\limits_{t=j+1}^{\hat{n}}a_{j,t}r_k^{(t)}\right|$.
It follows that 
$\theta_j<1$ and $\lim\limits_{k\rightarrow\infty} r_{k}^{(\hat{n})}=0$.
\end{proof}

\begin{rem}\label{rem:BB-convergence}
As $\hat{A}$ is positive definite, so is $\hat{A}^{-1}$. 
Let $\tilde{\lambda}_j=\tilde{u}_j+{\rm i}\tilde{v}_j~(1\le j\le \hat{n})$ be the eigenvalues of $\hat{A}^{-1}$. Clearly, $\frac{1}{\tilde{\lambda}_j}=\frac{\tilde{u}_j-{\rm i}\tilde{v}_j}{\tilde{u}_j^2+\tilde{v}_j^2}$ is an eigenvalue of $\hat{A}$. Then $\max\limits_{1\le j\le \hat{n}}\frac{|\lambda_j|^2}{u_j}
=\frac{1}{\min_{1\le j\le \hat{n}}\tilde{u}_j}$. 
This, along with
\begin{align*}
&\lambda_{\min}(W)
=2\lambda_{\min}\left((\hat{A}+\hat{A}^T)^{-1}\hat{A}^T\hat{A}\right)
=2\lambda_{\min}\left(\hat{A}(\hat{A}+\hat{A}^T)^{-1}\hat{A}^T\right)\\
&=2\lambda_{\min}\left((\hat{A}^{-1}+\hat{A}^{-T})^{-1}\right)
=\frac{2}{\lambda_{\max}(\hat{A}^{-1}+\hat{A}^{-T})},
\end{align*}
shows that condition \eqref{con1} is equivalent to 
$\lambda_{\max}(\hat{A}^{-1}+\hat{A}^{-T})<4\min\limits_{1\le j\le \hat{n}}\tilde{u}_j.$
Note that $\min\limits_{1\le j\le \hat{n}}\tilde{u}_j\ge\tfrac{1}{2}\lambda_{\min}(\hat{A}^{-1}+\hat{A}^{-T})$,%
\footnote{For any $j=1,\ldots,\hat{n}$, let $\tilde{x}_j$ be the eigenvector of $\hat{A}^{-1}$ corresponding to $\tilde{\lambda}_j$. Then we have $\tilde{\lambda}_j=\tfrac{\tilde{x}_j^*\hat{A}^{-1}\tilde{x}_j}{\tilde{x}_j^*\tilde{x}_j}$. Since $\hat{A}^{-1}+\hat{A}^{-T}$ is SPD, it gives
$
\tilde{u}_j=\frac12\left(\tilde{\lambda}_j+\tilde{\lambda}_j^*\right)
=\frac12\left(\frac{\tilde{x}_j^*\hat{A}^{-1}\tilde{x}_j}{\tilde{x}_j^*\tilde{x}_j}+\frac{\tilde{x}_j^*\hat{A}^{-T}\tilde{x}_j}{\tilde{x}_j^*\tilde{x}_j}\right)
=\frac{\tilde{x}_j^*\left(\hat{A}^{-1}+\hat{A}^{-T}\right)\tilde{x}_j}{2\tilde{x}_j^*\tilde{x}_j}
\ge\frac{1}{2}\lambda_{\min}(\hat{A}^{-1}+\hat{A}^{-T}).
$} 
so the above inequality can be reinforced as
$\lambda_{\max}(\hat{A}^{-1}+\hat{A}^{-T})<2\lambda_{\min}(\hat{A}^{-1}+\hat{A}^{-T}).$
When $\hat{A}$ is SPD, it reduces to $\lambda_{\max}(\hat{A})<2\lambda_{\min}(\hat{A})$, which is the same as the convergence condition of the preconditioned BB method for SPD linear systems \cite{molina1996}. This means that our condition \eqref{con1} is weaker than that of \cite{molina1996}. 
\end{rem}

When $G$ is UPD, so is $M$ in \eqref{a51}. Combining \Cref{thbb} with the convergence conditions of \Cref{inalg1} gives the following result.

\begin{theorem}\label{con_algbb}
Suppose $G\in \R^{n\times n}$ is UPD. For any SPD  $Q\in \R^{m\times m}$ and any $\omega>0$, let $M$ be defined by \eqref{a51} and $\lambda_j=u_j+{\rm i}v_j~(1\le j\le n+m)$ be its $n+m$ eigenvalues. If 
\begin{equation}\label{algbb_condition}
  \smash[t]{
    \max_{1\le j\le n+m}\frac{|\lambda_j|^2}{u_j}<\frac{4}{\lambda_{\max}\left(M^{-1}+M^{-T}\right)},}
\end{equation}
then for sufficiently small $\delta$, $\{x_{k},y_k\}$ produced by \Cref{kktbb1} converges to a solution of~\eqref{saddle1}.
\end{theorem}

\begin{rem}
    The residuals generated by the BB method, even for SPD linear systems, are strong nonmonotonic, which poses a challenge for the convergence \cite{raydan1997,dai2002r}. This is also the reason why the convergence of \Cref{kktbb1} is intricate. Our convergence analysis of \Cref{kktbb1} by ensuring a decrease of $\|r_{k}\|_*$ is quite stringent, relying on a rather strong assumption \eqref{algbb_condition}. The nonmonotonic behavior of $\|r_k\|$ in \Cref{fig:ns3-ns5,fig:ns6-sd1} also indicates that the choices of $\omega$ in our numerical experiments do not meet \eqref{algbb_condition}. Thus, there is significant room for improving the convergence of the BB method for UPD linear systems and \Cref{kktbb1}.
\end{rem}

\begin{rem}
    Although assumption \eqref{algbb_condition} is strong, it is still possible to choose $\omega$ to satisfy it. Indeed, consider the special case $n=m=1$ and $M=\pmat{ a & b \\ -b & \omega}$ with $a>0$ and $b\in\R$. Since 
    $
    M^{-1}=\frac{1}{a\omega+b^2}\pmat{ \omega & -b \\ b & a},
    $
    we have 
    $$
    \lambda_{\min}\left(M^{-1}+M^{-T}\right) 
    = \frac{2\min\{a, \omega\}}{a\omega+b^2}
    \quad\mbox{and}\quad
    \lambda_{\max}\left(M^{-1}+M^{-T}\right) 
    = \frac{2\max\{a, \omega\}}{a\omega+b^2}.
    $$
    It follows from \Cref{rem:BB-convergence} that \eqref{algbb_condition} can be reinforced as 
    $\lambda_{\max}\left(M^{-1}+M^{-T}\right)\le 2\lambda_{\min}\left(M^{-1}+M^{-T}\right)$, namely, 
    $\max\{a, \omega\} \le 2\min\{a, \omega\}$.
    This implies that \eqref{algbb_condition} holds when $\omega\in \left[a/2,\,2a\right]$. 
    
    For the general case, we can apply preconditioning techniques to \eqref{a11} such that $M$ is well-conditioned. Preconditioning techniques for $M$ have been widely studied; see \cite{Benzi2005} and the references therein.
\end{rem}

\section{Numerical experiments}\label{sec:num}

We present the results of numerical tests to examine the feasibility and effectiveness of {\spalbb}. All experiments were run using MATLAB R2022b on a PC with an Intel(R) Core(TM) i7-1260P CPU @~2.10GHz and 32GB of RAM. The initial guess is taken to be the zero vector, and the algorithms are terminated when the number of iterations exceeds $10^5$ or $ {\rm Res}:= \|r_k\|/ \|r_0\| \le 10^{-6}$. We report the number of outer iterations, the total number of iterations (for {\spalbb}, it includes the number of inner iterations), the CPU time in seconds, and the final value of the relative residual, denoted by ``Oiter”, ``Titer”, ``CPU” and ``Res”.

In {\spalbb}, we set $Q=I$, the stopping criterion~\eqref{a52} for inner iterations with $\delta=0.5$ and $2$-norm, and tried $\omega=10^{-i}$ with $i=1$, $2$, $3$, $4$, $5$, denoted {\spalbb}($\omega$). We compared our method with BICGSTAB and restarted GMRES. We tested two restart values: $20$ and $50$, denoted GMRES(20) and GMRES(50).

\begin{ex}\label{ex1}
The steady-state Navier-Stokes equations are
\begin{equation}\label{Navier-Stokes}
  -\nu\nabla^2 {\bm u}+{\bm u}\cdot\nabla {\bm u}+\nabla p={\bm h} 
  \quad {\rm and} \quad
  {\rm div}\,{\bm u}=0, 
  \quad {\bm z}=(x,y)\in\Omega,
\end{equation}
where $\Omega\subseteq \R^2$ is a bounded domain, the vector field ${\bm u}$ represents the velocity in $\Omega$, $p$ represents pressure, and $\nu>0$ is the kinematic viscosity. The test problem is a model of the flow in a square cavity $\Omega=(-1,1)\times(-1,1)$ with the lid moving from left to right. A Dirichlet no-flow (zero velocity) condition is applied on the side and bottom boundaries, and the nonzero horizontal velocity on the lid is $\{y=1;-1\le x\le 1\mid u_x=1-x^4\}$.
\end{ex}

Finite element discretization of \eqref{Navier-Stokes} results in 
system \eqref{saddle1} with $G=\nu G_1+G_2$. Here $G_1$ is SPD and consists of a set of uncouple discrete Laplace operators, corresponding to diffusion, and $G_2$ is a discrete convection operator and is unsymmetric. Evidently, $G$ becomes more unsymmetric as $\nu$ decreases. Various methods have been developed for solving~\eqref{Navier-Stokes}. However, the convergence rates of some approaches deteriorate as $\nu$ decreases \cite{elman1999preconditioning}. Thus, for~\eqref{Navier-Stokes}, we test three small viscosity values of $\nu$: $0.005,\,0.01,\,0.05$.

\begin{figure}[t]
	\subfloat[$\nu=0.005$]{
		\includegraphics[width=.5\linewidth]{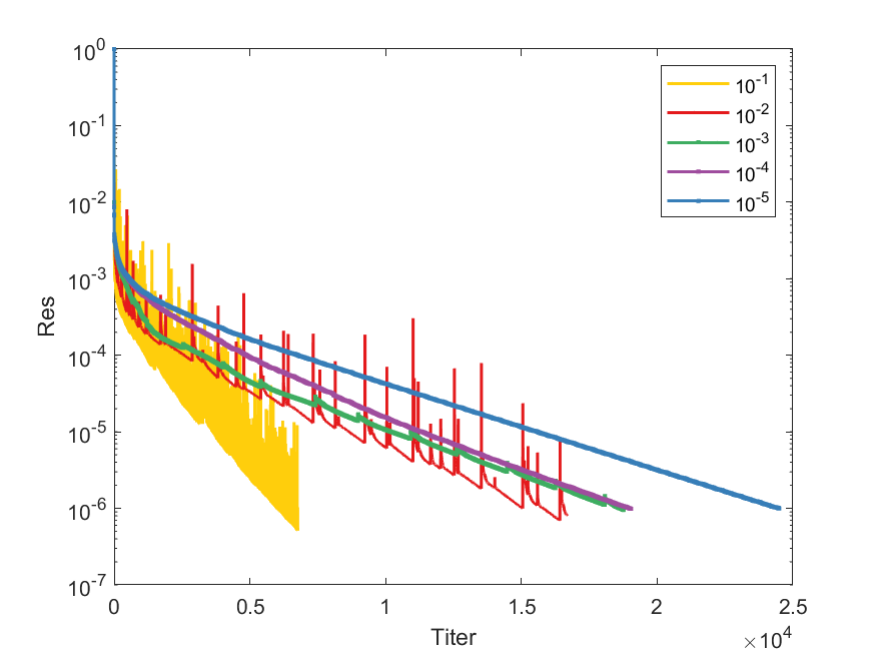}
        }
        \hfill
        \subfloat[$\nu=0.005$]{
		\includegraphics[width=.5\linewidth]{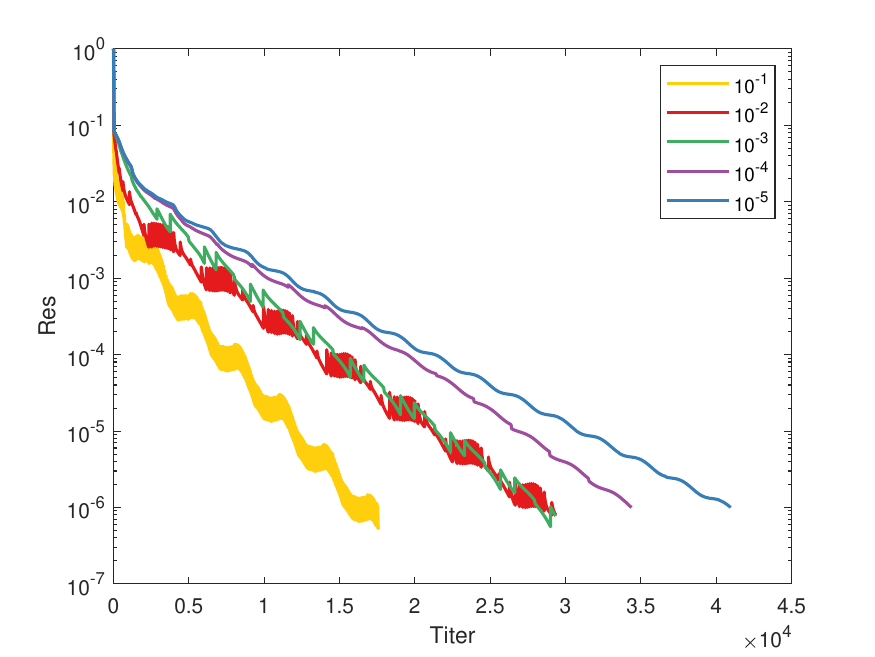}
        }
	\\
	\subfloat[$\nu=0.01$]{
		\includegraphics[width=.5\linewidth]{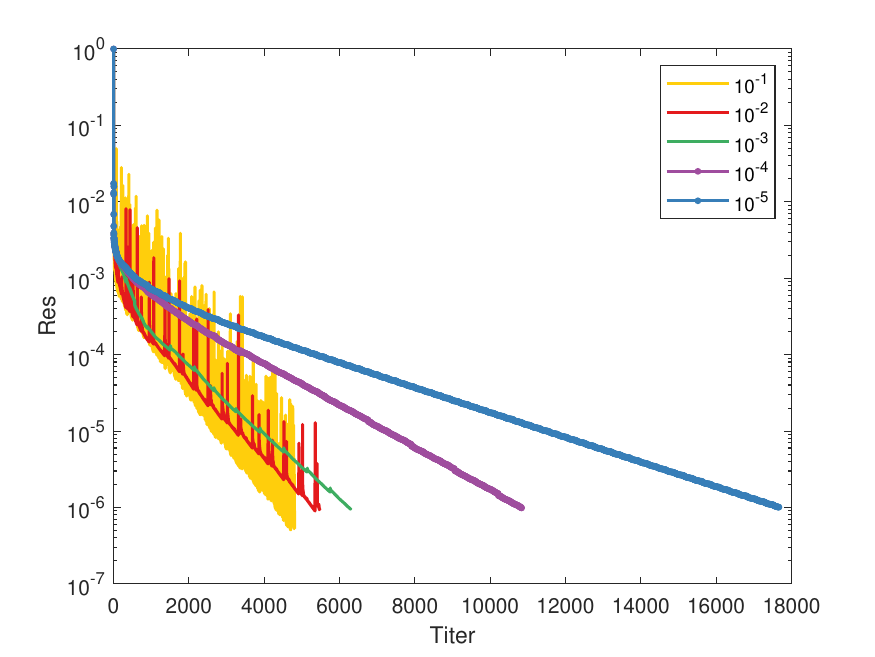}
        }
        \hfill
        \subfloat[$\nu=0.01$]{
		\includegraphics[width=.5\linewidth]{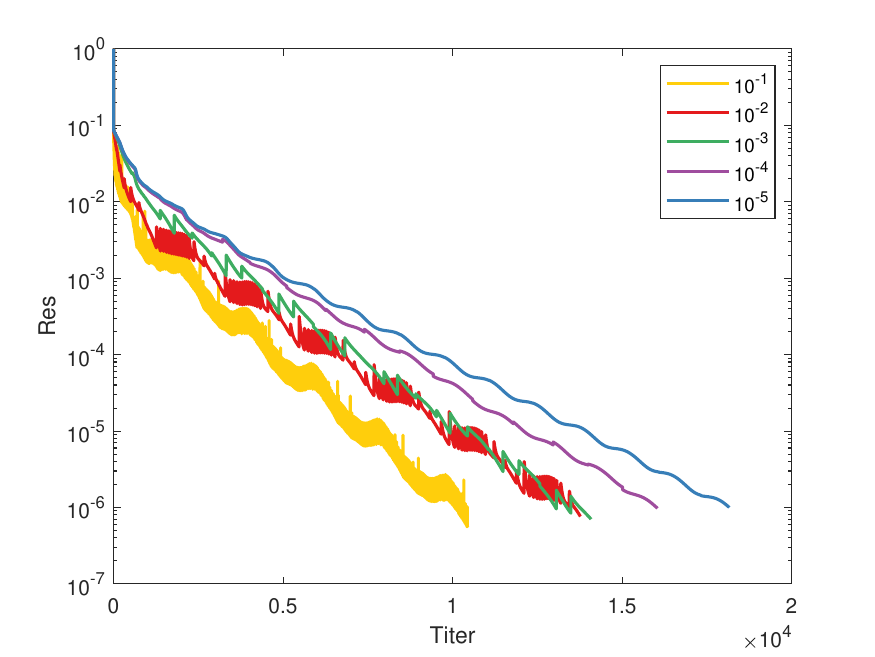}
        }
	\\
	\subfloat[$\nu=0.05$]{
		\includegraphics[width=.5\linewidth]{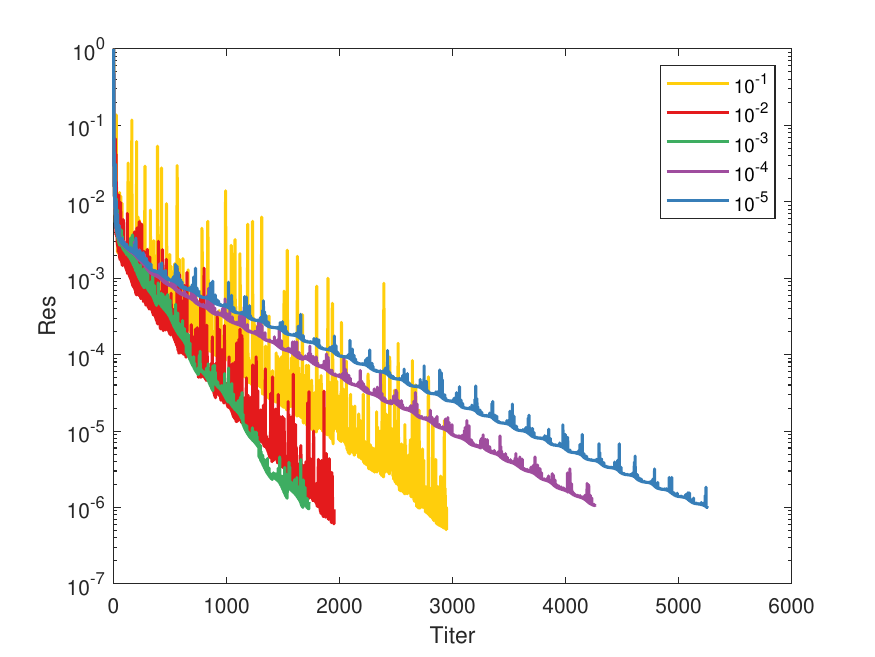}
	}
        \hfill
        \subfloat[$\nu=0.05$]{
		\includegraphics[width=.5\linewidth]{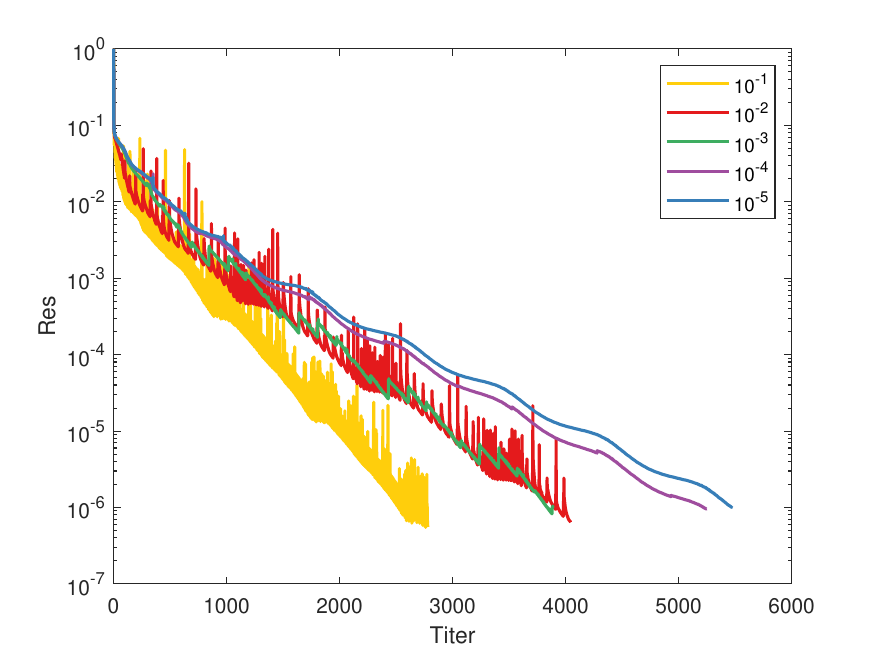}
	}
	\caption{Evolution of the relative residual of {\spalbb} tested on \Cref{ex1} (left) with $n=8450$, $m=1089$, and on \Cref{ex2} (right) with $n=8416$, $m=1096$ and $\omega$ as in \eqref{a11}.}
	\label{fig:ns3-ns5}
\end{figure}

\Cref{ex1} is a classical test problem used in fluid dynamics, known as driven-cavity flow. We discretize~\eqref{Navier-Stokes} using Picard iterations and the Q2--Q1 mixed finite element approximation \citep{elman2007} on uniform grids with grid parameter $h=2^{-6}$, $2^{-7}$, $2^{-8}$, $2^{-9}$. This discrete process can be accomplished by the IFISS software package \cite{elman2007,ifiss}. In this example, $G$ is UPD and $B$ is rank-deficient with rank $m-1$. Thus, the matrix in~\eqref{saddle1} is singular. The numerical results are reported in \Cref{tab1,tab2,tab3} and in the left-hand plots of \Cref{fig:ns3-ns5}, where ``-" means that the method failed to solve the problem and bold face indicates that the method performs best in terms of CPU time. It can be seen from \Cref{tab1,tab2,tab3} that the CPU time of all tested methods increases as $\nu$ decreases, and BICGSTAB and {\spalbb}(1) fail when $h=2^{-9}$ for $\nu=0.005$. The CPU time of {\spalbb} with $\omega\le10^{-2}$ is about 
half that of GMRES, and the best cases of {\spalbb} are only a third of GMRES for $h=2^{-9}$. The number of outer iterations of {\spalbb} decreases with $\omega$, which is consistent with~\Cref{rem:derome}. Nevertheless, the total number of iterations is not the least for $\omega=10^{-5}$.

\begin{table}[htbp]
\setlength\tabcolsep{4.1pt}
\centering\scriptsize
\caption{Numerical results for \Cref{ex1} with $\nu=0.005$. }\label{tab1}
\begin{tabular}{|c|r|r|r|r|r|r|r|r|r|r|r|} \hline
$h (n,m)$&\multicolumn{4}{c|}{$2^{-6}$ $(n=8,450,m=1,089)$}
&\multicolumn{4}{c|}{$2^{-7}$ $(n=33,282,m=4,225)$}\\\cline{2-9}
  &Oiter&Titer&CPU&RES&Oiter&Titer&CPU&RES\\\hline
BICGSTAB	&		&	6665.5	&	2.28 	&	$9.83$E$-07$	&		&	14211.5	&	22.75 	&	$8.95$E$-07$	\\
GMRES(20)	&		&	4221	&	2.55 	&	$9.99$E$-07$	&		&	7735	&	16.86 	&	$9.99$E$-07$	\\
GMRES(50)	&		&	4040	&	4.64 	&	$9.99$E$-07$	&		&	7162	&	27.94 	&	$1.00$E$-06$	\\
{\spalbb}(1)	&	1594	&	6765	&	\pmb{1.23} 	&	$9.97$E$-07$	&	22057	&	50212	&	42.16 	&	$9.93$E$-07$	\\
{\spalbb}(2)	&	51	&	16705	&	2.83 	&	$9.98$E$-07$	&	243	&	18537	&	\pmb{15.42} 	&	$9.83$E$-07$	\\
{\spalbb}(3)	&	17	&	18762	&	3.18 	&	$1.00$E$-06$	&	27	&	24084	&	20.45 	&	$7.70$E$-07$	\\
{\spalbb}(4)	&	14	&	19036	&	3.22 	&	$9.99$E$-07$	&	15	&	22801	&	18.78 	&	$1.00$E$-06$	\\
{\spalbb}(5)	&	14	&	24496	&	4.28 	&	$1.00$E$-06$	&	14	&	35119	&	29.70 	&	$1.00$E$-06$	\\
\hline

$h (n,m)$&\multicolumn{4}{c|}{$2^{-8}$ $(n=132,098,m=16,641)$}
&\multicolumn{4}{c|}{$2^{-9}$ $(n=526,338,m=66,049)$}\\\cline{2-9}
  &Oiter&Titer&CPU&RES&Oiter&Titer&CPU&RES\\\hline
BICGSTAB	&		&	39440.5	&	378.63 	&	$8.31$E$-07$	&		&	-	&	-	&	-	\\
GMRES(20)	&		&	15638	&	138.37 	&	$1.00$E$-06$	&		&	37240	&	1563.36 	&	$1.00$E$-06$	\\
GMRES(50)	&		&	13265	&	179.20 	&	$1.00$E$-06$	&		&	25858	&	1518.14 	&	$1.00$E$-06$	\\
{\spalbb}(1)	&	55401	&	77589	&	374.99 	&	$9.99$E$-07$	&	-	&	-	&	-	&	-	\\
{\spalbb}(2)	&	2412	&	17982	&	\pmb{85.91} 	&	$1.00$E$-06$	&	14384	&	38247	&	808.82 	&	$9.94$E$-07$	\\
{\spalbb}(3)	&	51	&	26095	&	127.44 	&	$1.00$E$-06$	&	340	&	26085	&	\pmb{539.99} 	&	$9.69$E$-07$	\\
{\spalbb}(4)	&	17	&	28805	&	138.84 	&	$1.00$E$-06$	&	23	&	36728	&	778.87 	&	$1.00$E$-06$	\\
{\spalbb}(5)	&	14	&	32707	&	165.62 	&	$1.00$E$-06$	&	14	&	36635	&	783.08 	&	$1.00$E$-06$	\\
\hline
\end{tabular}
\end{table}

\begin{table}[htbp]
\setlength\tabcolsep{4.4pt}
\centering\scriptsize
\caption{Numerical results for \Cref{ex1} with $\nu=0.01$. }\label{tab2}
\begin{tabular}{|c|r|r|r|r|r|r|r|r|r|r|r|} \hline
$h (n,m)$&\multicolumn{4}{c|}{$2^{-6}$ $(n=8,450,m=1,089)$}
&\multicolumn{4}{c|}{$2^{-7}$ $(n=33,282,m=4,225)$}\\\cline{2-9}
  &Oiter&Titer&CPU&RES&Oiter&Titer&CPU&RES\\\hline
BICGSTAB	&		&	-	&	-	&	-	&		&	7211.5	&	11.63 	&	$9.90$E$-07$	\\
GMRES(20)	&		&	2453	&	1.60 	&	$9.98$E$-07$	&		&	5127	&	10.15 	&	$9.99$E$-07$	\\
GMRES(50)	&		&	2237	&	2.63 	&	$9.98$E$-07$	&		&	4422	&	14.82 	&	$1.00$E$-06$	\\
{\spalbb}(1)	&	1649	&	4803	&	\pmb{0.87} 	&	$9.96$E$-07$	&	8611	&	17952	&	15.60 	&	$9.49$E$-07$	\\
{\spalbb}(2)	&	38	&	5475	&	0.91 	&	$1.00$E$-06$	&	411	&	9982	&	8.26 	&	$1.00$E$-06$	\\
{\spalbb}(3)	&	16	&	6295	&	0.98 	&	$1.00$E$-06$	&	22	&	8554	&	\pmb{6.94} 	&	$1.00$E$-06$	\\
{\spalbb}(4)	&	15	&	10835	&	1.76 	&	$1.00$E$-06$	&	15	&	10388	&	8.35 	&	$1.00$E$-06$	\\
{\spalbb}(5)	&	15	&	17674	&	2.89 	&	$1.00$E$-06$	&	15	&	24493	&	20.09 	&	$1.00$E$-06$	\\\hline

$h (n,m)$&\multicolumn{4}{c|}{$2^{-8}$ $(n=132,098,m=16,641)$}
&\multicolumn{4}{c|}{$2^{-9}$ $(n=526,338,m=66,049)$}\\\cline{2-9}
  &Oiter&Titer&CPU&RES&Oiter&Titer&CPU&RES\\\hline
BICGSTAB	&		&	17765.5	&	161.00 	&	$8.20$E$-07$	&		&	29320	&	1024.60 	&	$1.46$E$-05$	\\
GMRES(20)	&		&	13579	&	113.90 	&	$1.00$E$-06$	&		&	34224	&	1186.20 	&	$1.00$E$-06$	\\
GMRES(50)	&		&	9419	&	118.86 	&	$1.00$E$-06$	&		&	22827	&	1140.94 	&	$1.00$E$-06$	\\
{\spalbb}(1)	&	25643	&	36571	&	165.72 	&	$9.99$E$-07$	&	77225	&	87388	&	1632.85 	&	$9.94$E$-07$	\\
{\spalbb}(2)	&	2154	&	14993	&	67.21 	&	$1.00$E$-06$	&	8277	&	31640	&	581.72 	&	$9.99$E$-07$	\\
{\spalbb}(3)	&	42	&	13434	&	\pmb{60.32} 	&	$9.99$E$-07$	&	554	&	22207	&	402.99 	&	$1.00$E$-06$	\\
{\spalbb}(4)	&	16	&	14113	&	65.49 	&	$9.99$E$-07$	&	22	&	21902	&	395.92 	&	$9.14$E$-07$	\\
{\spalbb}(5)	&	15	&	21503	&	101.01 	&	$1.00$E$-06$	&	15	&	20219	&	\pmb{369.43} 	&	$1.00$E$-06$	\\\hline

\end{tabular}
\end{table}

\begin{table}[htbp]
\setlength\tabcolsep{4.6pt}
\centering\scriptsize
\caption{Numerical results for \Cref{ex1} with $\nu=0.05$. }\label{tab3}
\begin{tabular}{|c|r|r|r|r|r|r|r|r|r|r|r|} \hline
$h (n,m)$&\multicolumn{4}{c|}{$2^{-6}$ $(n=8,450,m=1,089)$}
&\multicolumn{4}{c|}{$2^{-7}$ $(n=33,282,m=4,225)$}\\\cline{2-9}
  &Oiter&Titer&CPU&RES&Oiter&Titer&CPU&RES\\\hline
BICGSTAB	&		&	917.5	&	0.33 	&	$9.32$E$-07$	&		&	1648.5	&	2.72 	&	$9.61$E$-07$	\\
GMRES(20)	&		&	1508	&	0.94 	&	$9.96$E$-07$	&		&	4507	&	9.58 	&	$9.98$E$-07$	\\
GMRES(50)	&		&	989	&	1.15 	&	$9.96$E$-07$	&		&	2769	&	11.11 	&	$1.00$E$-06$	\\
{\spalbb}(1)	&	782	&	2947	&	0.54 	&	$9.96$E$-07$	&	3142	&	9527	&	7.89 	&	$1.00$E$-06$	\\
{\spalbb}(2)	&	79	&	1952	&	0.34 	&	$9.95$E$-07$	&	300	&	4762	&	3.84 	&	$9.98$E$-07$	\\
{\spalbb}(3)	&	19	&	1729	&	\pmb{0.28} 	&	$9.97$E$-07$	&	36	&	3272	&	\pmb{2.69} 	&	$6.68$E$-07$	\\
{\spalbb}(4)	&	17	&	4261	&	0.67 	&	$9.98$E$-07$	&	15	&	3878	&	3.21 	&	$9.95$E$-07$	\\
{\spalbb}(5)	&	17	&	5252	&	0.86 	&	$1.00$E$-06$	&	16	&	8563	&	7.24 	&	$9.97$E$-07$	\\
\hline

$h (n,m)$&\multicolumn{4}{c|}{$2^{-8}$ $(n=132,098,m=16,641)$}
&\multicolumn{4}{c|}{$2^{-9}$ $(n=526,338,m=66,049)$}\\\cline{2-9}
  &Oiter&Titer&CPU&RES&Oiter&Titer&CPU&RES\\\hline
BICGSTAB	&		&	3253.5	&	30.99 	&	$8.58$E$-07$	&		&	6277.5	&	252.29 	&	$9.54$E$-07$	\\
GMRES(20)	&		&	11071	&	99.96 	&	$1.00$E$-06$	&		&	21205	&	797.99 	&	$1.00$E$-06$	\\
GMRES(50)	&		&	7775	&	104.78 	&	$1.00$E$-06$	&		&	18460	&	1048.97 	&	$1.00$E$-06$	\\
{\spalbb}(1)	&	10731	&	30596	&	141.61 	&	$9.99$E$-07$	&	-	&	-	&	-	&	-	\\
{\spalbb}(2)	&	1005	&	13335	&	62.92 	&	$1.00$E$-06$	&	3330	&	34174	&	700.69 	&	$9.99$E$-07$	\\
{\spalbb}(3)	&	100	&	7869	&	36.87 	&	$9.99$E$-07$	&	339	&	21075	&	433.98 	&	$1.00$E$-06$	\\
{\spalbb}(4)	&	19	&	5156	&	\pmb{24.64} 	&	$7.36$E$-07$	&	37	&	10542	&	\pmb{216.78} 	&	$9.99$E$-07$	\\
{\spalbb}(5)	&	17	&	12476	&	60.19 	&	$9.99$E$-07$	&	16	&	11154	&	231.77 	&	$9.99$E$-07$	\\
\hline
\end{tabular}
\end{table}


\begin{ex}\label{ex2}
We consider the steady-state Navier-Stokes equations~\eqref{Navier-Stokes}, where the domain $\Omega$ is a rectangular region $(0,8)\times(-1,1)$ generated by deleting the square $(7/4,9/4)\times
(-1/4,1/4)$. This test problem is a model of the flow in a rectangular channel with an obstacle. A Poiseuille profile is imposed on the inflow boundary $\{x=0; -1\le y\le 1\}$, and a Dirichlet no-flow condition is imposed on the obstruction and on the top and bottom walls. A Neumann condition is applied at the outflow boundary that automatically sets the mean outflow pressure to zero.
\end{ex}

In our tests, we set $\nu=0.005,\,0.01,\,0.05$ and discretize the Navier-Stokes equations~\eqref{Navier-Stokes} using Picard iterations and the Q2--Q1 mixed finite element approximation \citep{elman2007} on uniform grids with grid parameter $h=2^{-5},\,2^{-6},\,2^{-7},\,2^{-8}$. This discretization was accomplished using IFISS \citep{elman2007,ifiss}. The resulting matrices have $G$ UPD and $B$ full column rank. The numerical results are reported in \Cref{tab4,tab5,tab6} and \Cref{fig:ns3-ns5}. \Cref{tab4,tab5,tab6} show that all choices of $\omega$ are successful in solving the tested problems, and, in terms of CPU time, $\omega=10^{-1}$ and $10^{-2}$ perform better than other choices. Although BICGSTAB requires the least CPU time for $\nu=0.05$, it fails for $\nu=0.005$ and $\nu=0.01$ with $h=2^{-5},\,2^{-8}$. The CPU time for every {\spalbb} test is less than for GMRES, and the best case of {\spalbb} takes about half the time of GMRES. Overall, {\spalbb} is more stable and efficient.


\begin{table}[htbp]
\setlength\tabcolsep{4.6pt}
\centering\scriptsize
\caption{Numerical results for \Cref{ex2} with $\nu=0.005$. }\label{tab4}
\begin{tabular}{|c|r|r|r|r|r|r|r|r|r|r|r|} \hline
$h (n,m)$&\multicolumn{4}{c|}{$2^{-5}$ $(n=8,416,m=1,096)$}
&\multicolumn{4}{c|}{$2^{-6}$ $(n=32,960,m=4,208)$}\\\cline{2-9}
  &Oiter&Titer&CPU&RES&Oiter&Titer&CPU&RES\\\hline
BICGSTAB	&		&	-	&	-	&	-	&		&	-	&	-	&	-	\\
GMRES(20)	&		&	6319	&	3.60 	&	$9.99$E$-07$	&		&	11462	&	25.25 	&	$9.99$E$-07$	\\
GMRES(50)	&		&	5540	&	6.10 	&	$9.98$E$-07$	&		&	11386	&	42.46 	&	$1.00$E$-06$	\\
{\spalbb}(1)	&	632	&	17567	&	\pmb{2.78} 	&	$1.00$E$-06$	&	2211	&	16590	&	\pmb{12.88} 	&	$1.00$E$-06$	\\
{\spalbb}(2)	&	110	&	29351	&	4.95	&	$1.00$E$-06$	&	330	&	34576	&	27.70 	&	$9.99$E$-07$	\\
{\spalbb}(3)	&	30	&	29332	&	4.64 	&	$1.00$E$-06$	&	61	&	34000	&	26.83 	&	$8.13$E$-07$	\\
{\spalbb}(4)	&	18	&	34374	&	5.25 	&	$1.00$E$-06$	&	20	&	34842	&	27.95 	&	$9.99$E$-07$	\\
{\spalbb}(5)	&	18	&	40951	&	6.20 	&	$1.00$E$-06$	&	17	&	44791	&	36.40 	&	$1.00$E$-06$	\\
\hline

$h (n,m)$&\multicolumn{4}{c|}{$2^{-7}$ $(n=130,432,m=16,480)$}
&\multicolumn{4}{c|}{$2^{-8}$ $(n=518,912,m=65,216)$}\\\cline{2-9}
  &Oiter&Titer&CPU&RES&Oiter&Titer&CPU&RES\\\hline
BICGSTAB	&		&	-	&	-	&	-	&		&	-	&	-	&	-	\\
GMRES(20)	&		&	20442	&	173.54 	&	$1.00$E$-06$	&		&	39863	&	1385.75 	&	$1.00$E$-06$	\\
GMRES(50)	&		&	20511	&	276.89 	&	$1.00$E$-06$	&		&	38382	&	2019.57 	&	$1.00$E$-06$	\\
{\spalbb}(1)	&	9013	&	23219	&	\pmb{102.78} 	&	$1.00$E$-06$	&	43791	&	61073	&	1155.48 	&	$1.00$E$-06$	\\
{\spalbb}(2)	&	917	&	46379	&	211.82 	&	$9.97$E$-07$	&	2765	&	35310	&	\pmb{651.49} 	&	$1.00$E$-06$	\\
{\spalbb}(3)	&	135	&	44805	&	214.76 	&	$1.00$E$-06$	&	361	&	60088	&	1099.24 	&	$1.00$E$-06$	\\
{\spalbb}(4)	&	32	&	46194	&	221.46 	&	$9.10$E$-07$	&	65	&	59081	&	1080.07 	&	$8.08$E$-07$	\\
{\spalbb}(5)	&	16	&	52684	&	255.82 	&	$1.00$E$-06$	&	19	&	57737	&	1064.29 	&	$1.00$E$-06$	\\
\hline
\end{tabular}
\end{table}

\begin{table}[htbp]
\setlength\tabcolsep{4.6pt}
\centering\scriptsize
\caption{Numerical results for \Cref{ex2} with $\nu=0.01$. }\label{tab5}
\begin{tabular}{|c|r|r|r|r|r|r|r|r|r|r|r|} \hline
$h (n,m)$&\multicolumn{4}{c|}{$2^{-5}$ $(n=8,416,m=1,096)$}
&\multicolumn{4}{c|}{$2^{-6}$ $(n=32,960,m=4,208)$}\\\cline{2-9}
  &Oiter&Titer&CPU&RES&Oiter&Titer&CPU&RES\\\hline
BICGSTAB	&		&	-	&	-	&	-	&		&	5336	&	\pmb{7.87} 	&	$9.89$E$-07$	\\
GMRES(20)	&		&	5145	&	2.98 	&	$9.99$E$-07$	&		&	9162	&	18.05 	&	$1.00$E$-06$	\\
GMRES(50)	&		&	4904	&	5.59 	&	$1.00$E$-06$	&		&	9446	&	37.68 	&	$1.00$E$-06$	\\
{\spalbb}(1)	&	604	&	10445	&	\pmb{1.54} 	&	$9.99$E$-07$	&	2455	&	10813	&	8.26 	&	$9.99$E$-07$	\\
{\spalbb}(2)	&	108	&	13769	&	2.06 	&	$9.42$E$-07$	&	294	&	20425	&	15.79 	&	$9.77$E$-07$	\\
{\spalbb}(3)	&	30	&	14084	&	2.95 	&	$6.63$E$-07$	&	57	&	19994	&	14.99 	&	$9.99$E$-07$	\\
{\spalbb}(4)	&	18	&	16042	&	2.50 	&	$9.98$E$-07$	&	20	&	21636	&	16.45 	&	$7.17$E$-07$	\\
{\spalbb}(5)	&	18	&	18160	&	2.75 	&	$1.00$E$-06$	&	17	&	26448	&	20.12 	&	$1.00$E$-06$	\\
\hline

$h (n,m)$&\multicolumn{4}{c|}{$2^{-7}$ $(n=130,432,m=16,480)$}
&\multicolumn{4}{c|}{$2^{-8}$ $(n=518,912,m=65,216)$}\\\cline{2-9}
  &Oiter&Titer&CPU&RES&Oiter&Titer&CPU&RES\\\hline
BICGSTAB	&		&	11686.5	&	104.74 	&	$7.60$E$-07$	&		&	-	&	-	&	-	\\
GMRES(20)	&		&	17521	&	151.67 	&	$1.00$E$-06$	&		&	34452	&	1190.29 	&	$1.00$E$-06$	\\
GMRES(50)	&		&	17430	&	237.68 	&	$1.00$E$-06$	&		&	33304	&	1629.56 	&	$1.00$E$-06$	\\
{\spalbb}(1)	&	8001	&	17667	&	\pmb{79.71} 	&	$9.98$E$-07$	&	30579	&	43828	&	851.71 	&	$9.99$E$-07$	\\
{\spalbb}(2)	&	818	&	26085	&	115.24 	&	$1.00$E$-06$	&	2576	&	28375	&	\pmb{496.34} 	&	$1.00$E$-06$	\\
{\spalbb}(3)	&	138	&	29649	&	143.63 	&	$9.98$E$-07$	&	308	&	40365	&	719.53 	&	$9.99$E$-07$	\\
{\spalbb}(4)	&	31	&	28824	&	133.78 	&	$6.90$E$-07$	&	67	&	47090	&	870.44 	&	$9.64$E$-07$	\\
{\spalbb}(5)	&	17	&	40333	&	172.89 	&	$1.00$E$-06$	&	20	&	48639	&	909.21 	&	$7.29$E$-07$	\\
\hline
\end{tabular}
\end{table}

\begin{table}[htbp]
\setlength\tabcolsep{4.6pt}
\centering\scriptsize
\caption{Numerical results for \Cref{ex2} with $\nu=0.05$. }\label{tab6}
\begin{tabular}{|c|r|r|r|r|r|r|r|r|r|r|r|} \hline
$h (n,m)$&\multicolumn{4}{c|}{$2^{-5}$ $(n=8,416,m=1,096)$}
&\multicolumn{4}{c|}{$2^{-6}$ $(n=32,960,m=4,208)$}\\\cline{2-9}
  &Oiter&Titer&CPU&RES&Oiter&Titer&CPU&RES\\\hline
BICGSTAB	&		&	1127.5	&	\pmb{0.39} 	&	$7.52$E$-07$	&		&	2187.5	&	\pmb{3.24} 	&	$5.81$E$-07$	\\
GMRES(20)	&		&	2420	&	2.63 	&	$9.97$E$-07$	&		&	5288	&	10.27 	&	$1.00$E$-06$	\\
GMRES(50)	&		&	2686	&	4.79 	&	$9.99$E$-07$	&		&	5168	&	17.39 	&	$9.99$E$-07$	\\
{\spalbb}(1)	&	616	&	2780	&	0.47 	&	$9.92$E$-07$	&	1877	&	6454	&	4.69 	&	$9.99$E$-07$	\\
{\spalbb}(2)	&	90	&	4043	&	0.59 	&	$9.98$E$-07$	&	222	&	6876	&	4.96 	&	$9.96$E$-07$	\\
{\spalbb}(3)	&	28	&	3888	&	0.63 	&	$1.00$E$-06$	&	47	&	6949	&	5.06 	&	$9.98$E$-07$	\\
{\spalbb}(4)	&	19	&	5246	&	0.79 	&	$9.97$E$-07$	&	21	&	8136	&	6.06 	&	$1.00$E$-06$	\\
{\spalbb}(5)	&	18	&	5470	&	0.83 	&	$9.97$E$-07$	&	18	&	9266	&	6.89 	&	$9.99$E$-07$	\\
\hline

$h (n,m)$&\multicolumn{4}{c|}{$2^{-7}$ $(n=130,432,m=16,480)$}
&\multicolumn{4}{c|}{$2^{-8}$ $(n=518,912,m=65,216)$}\\\cline{2-9}
  &Oiter&Titer&CPU&RES&Oiter&Titer&CPU&RES\\\hline
BICGSTAB	&		&	4516.5	&	\pmb{39.77} 	&	$9.48$E$-07$	&		&	9366.5	&	\pmb{340.85} 	&	$8.40$E$-07$	\\
GMRES(20)	&		&	10778	&	89.17 	&	$9.99$E$-07$	&		&	20466	&	711.25 	&	$1.00$E$-06$	\\
GMRES(50)	&		&	10134	&	130.59 	&	$1.00$E$-06$	&		&	18845	&	921.70 	&	$9.99$E$-07$	\\
{\spalbb}(1)	&	6815	&	19080	&	83.73 	&	$9.98$E$-07$	&	22561	&	57886	&	1903.39 	&	$1.00$E$-06$	\\
{\spalbb}(2)	&	738	&	11894	&	53.23 	&	$9.96$E$-07$	&	2657	&	29081	&	536.54 	&	$9.99$E$-07$	\\
{\spalbb}(3)	&	100	&	12759	&	57.35 	&	$9.98$E$-07$	&	304	&	24471	&	451.52 	&	$1.00$E$-06$	\\
{\spalbb}(4)	&	28	&	13047	&	59.88 	&	$1.00$E$-06$	&	48	&	23291	&	434.16 	&	$1.00$E$-06$	\\
{\spalbb}(5)	&	18	&	15504	&	71.03 	&	$1.00$E$-06$	&	21	&	25483	&	469.27 	&	$9.98$E$-07$	\\
\hline
\end{tabular}
\end{table}


\begin{ex}\label{ex3}
We consider the steady-state Navier-Stokes equations~\eqref{Navier-Stokes}, where the domain $\Omega$ is a rectangular region $(-1,5)\times(-1,1)$ generated by deleting $(-1,0)\times(-1,-1/2)\cup(-1,0)\times(1/2,1)$. This test problem is a model of the flow in a symmetric step channel. A Poiseuille flow profile is imposed on the inflow boundary $\{x=-1; -1/2\le y\le1/2\}$, and a Dirichlet no-flow condition is imposed on the top and bottom walls and the boundaries of deleted parts. A Neumann condition is applied at the outflow boundary that sets the mean outflow pressure to zero.
\end{ex}

\begin{figure}[t]
	\subfloat[$\nu=0.005$]{
		\includegraphics[width=.5\linewidth]{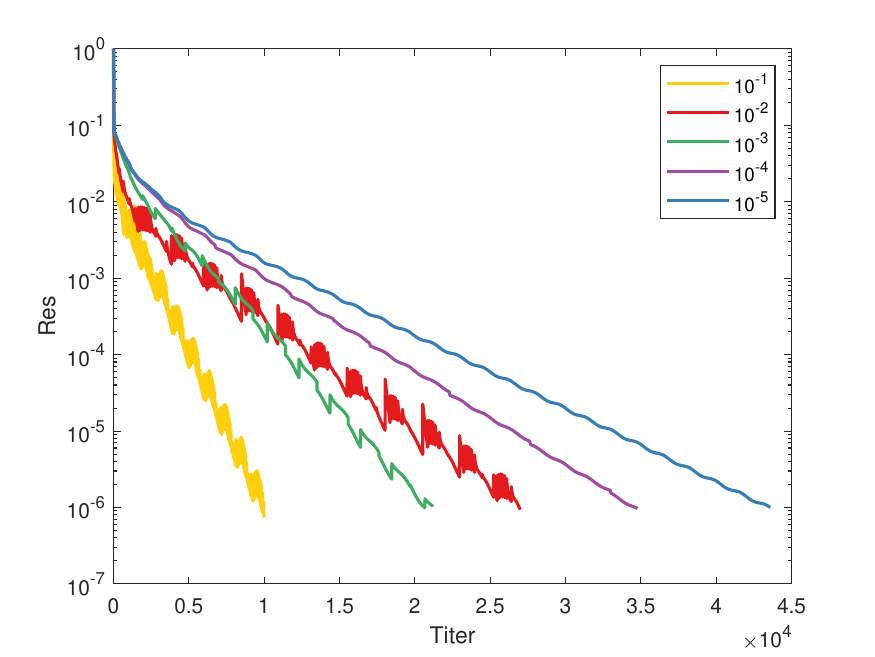}
	}
        \hfill
        \subfloat[$\nu=0.005$]{
		\includegraphics[width=.5\linewidth]{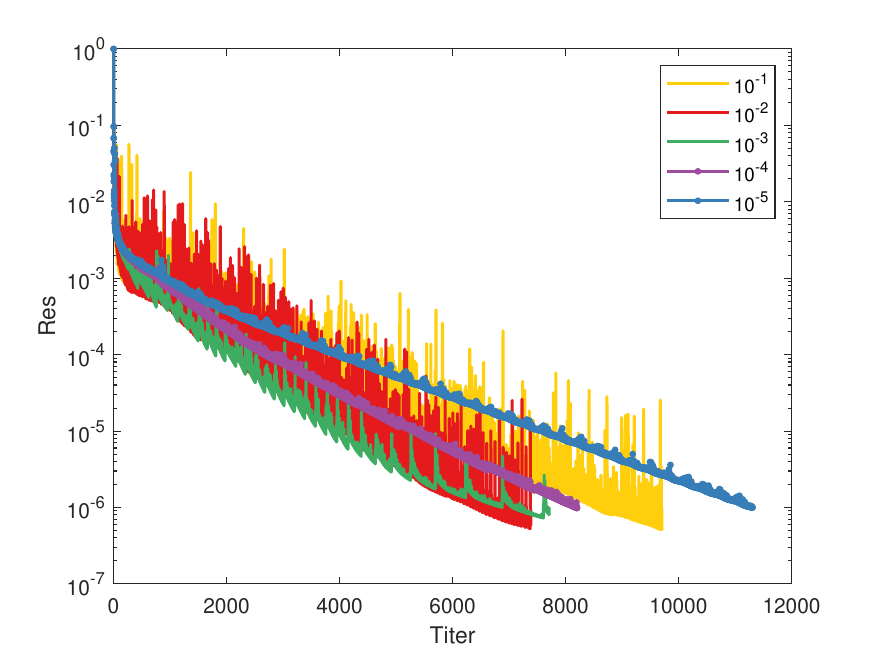}
	}
        \\
	\subfloat[$\nu=0.01$]{
		\includegraphics[width=.5\linewidth]{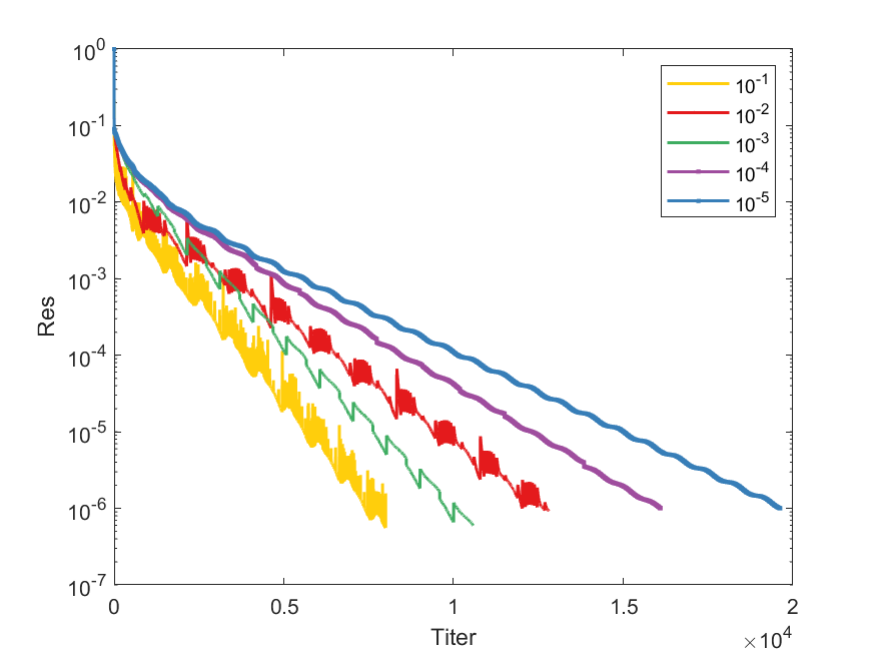}
	}
        \hfill
        \subfloat[$\nu=0.01$]{
		\includegraphics[width=.5\linewidth]{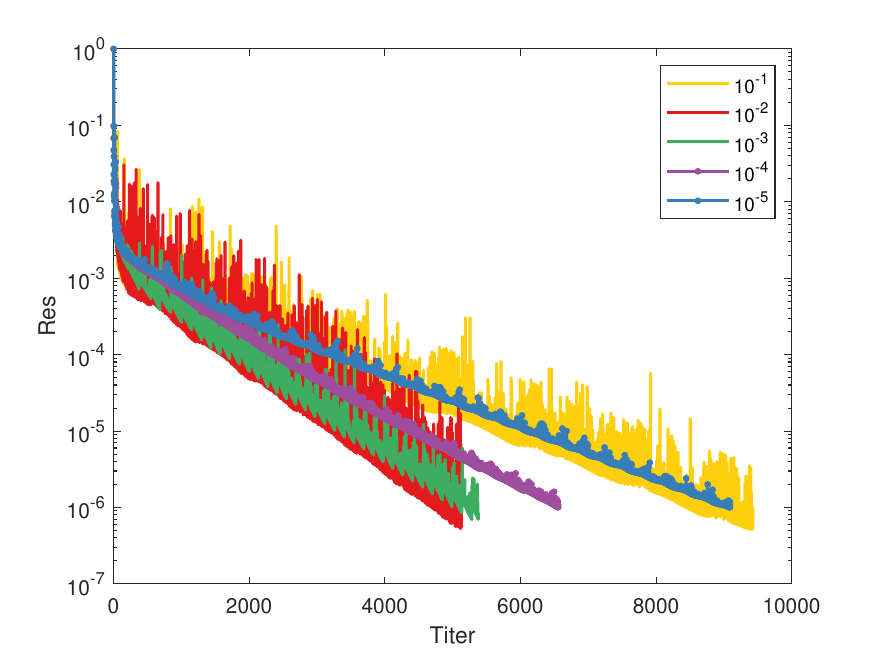}
        }
        \\
	\subfloat[$\nu=0.05$]{
		\includegraphics[width=.5\linewidth]{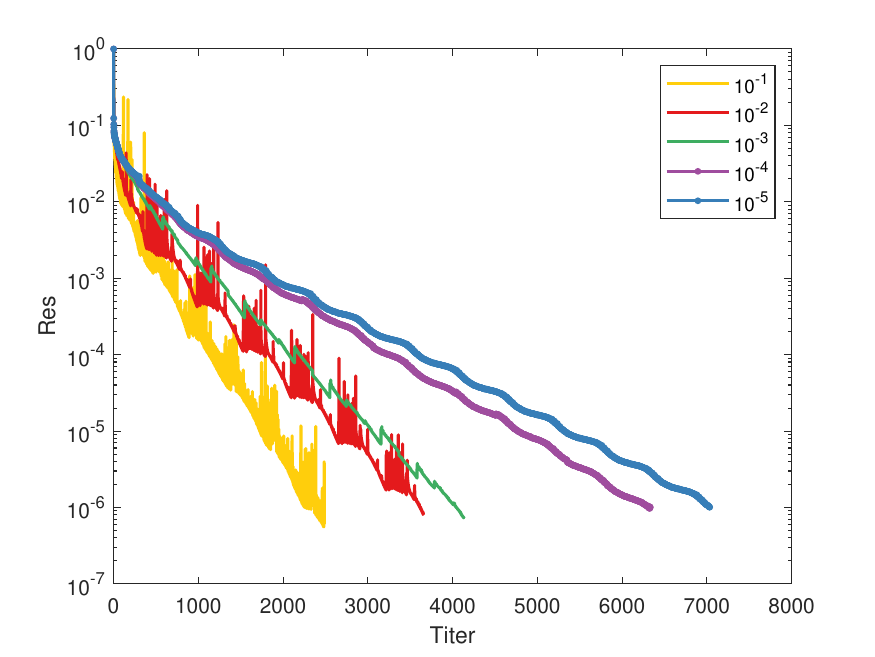}
	}
        \hfill
        \subfloat[$\nu=0.05$]{
		\includegraphics[width=.5\linewidth]{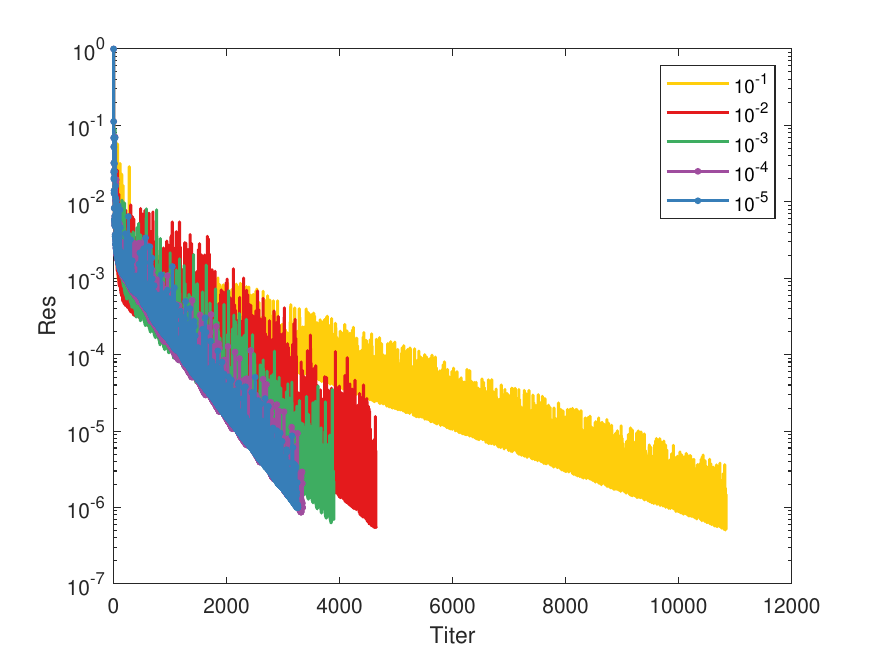}
	}
	\caption{Evolution of the relative residual of {\spalbb} tested on \Cref{ex3} (left) with $n=5890$, $m=769$, and on \Cref{ex:sd} with $n=12675$, $m=1089$ and $\omega$ as in \eqref{a11}.}
	\label{fig:ns6-sd1}
\end{figure}

The discretization of the Navier-Stokes equations~\eqref{Navier-Stokes} is done as in \Cref{ex2} with the same setting. In this example, $G$ is UPD and $B$ has full column rank. The numerical results are reported in \Cref{tab7,tab8,tab9} and in the left-hand plots of \Cref{fig:ns6-sd1}. As in \Cref{ex2}, all choices of $\omega$ solve the problems successfully, and BICGSTAB performs best in the case of $\nu=0.05$. Except for $\nu=0.05$ and $\nu=0.01$ with $h=2^{-6}$, {\spalbb} requires the least CPU time. Hence, \Cref{tab7,tab8,tab9} still demonstrate the efficiency of {\spalbb}.

\begin{table}[htbp]
\setlength\tabcolsep{4.6pt}
\centering\scriptsize
\caption{Numerical results for \Cref{ex3} with $\nu=0.005$. }\label{tab7}
\begin{tabular}{|c|r|r|r|r|r|r|r|r|r|r|r|} \hline
$h (n,m)$&\multicolumn{4}{c|}{$2^{-5}$ $(n=5,890,m=769)$}
&\multicolumn{4}{c|}{$2^{-6}$ $(n=23,042,m=2,945)$}\\\cline{2-9}
  &Oiter&Titer&CPU&RES&Oiter&Titer&CPU&RES\\\hline
BICGSTAB	&		&	-	&	-	&	-	&		&	-	&	-	&	-	\\
GMRES(20)	&		&	5662	&	2.23 	&	$9.98$E$-07$	&		&	12340	&	17.73 	&	$1.00$E$-06$	\\
GMRES(50)	&		&	3918	&	3.01 	&	$1.00$E$-06$	&		&	10459	&	26.94 	&	$9.99$E$-07$	\\
{\spalbb}(1)	&	510	&	10007	&	\pmb{0.96} 	&	$9.96$E$-07$	&	2392	&	13447	&	\pmb{6.73} 	&	$9.99$E$-07$	\\
{\spalbb}(2)	&	99	&	26993	&	2.44 	&	$9.44$E$-07$	&	322	&	36764	&	17.71 	&	$1.00$E$-06$	\\
{\spalbb}(3)	&	23	&	21194	&	1.92 	&	$1.00$E$-06$	&	54	&	29542	&	14.58 	&	$9.39$E$-07$	\\
{\spalbb}(4)	&	17	&	34776	&	3.20 	&	$1.00$E$-06$	&	20	&	36454	&	18.04 	&	$1.00$E$-06$	\\
{\spalbb}(5)	&	18	&	43570	&	3.94 	&	$1.00$E$-06$	&	17	&	47559	&	23.47 	&	$1.00$E$-06$	\\\hline

$h (n,m)$&\multicolumn{4}{c|}{$2^{-7}$ $(n=91,138,m=11,521)$}
&\multicolumn{4}{c|}{$2^{-8}$ $(n=362,498,m=45,569)$}\\\cline{2-9}
  &Oiter&Titer&CPU&RES&Oiter&Titer&CPU&RES\\\hline
BICGSTAB	&		&	-	&	-	&	-	&		&	-	&	-	&	-	\\
GMRES(20)	&		&	21340	&	119.38 	&	$1.00$E$-06$	&		&	41829	&	954.23 	&	$1.00$E$-06$	\\
GMRES(50)	&		&	22164	&	196.79 	&	$9.99$E$-07$	&		&	41230	&	1409.32 	&	$1.00$E$-06$	\\
{\spalbb}(1)	&	9838	&	20630	&	\pmb{50.41} 	&	$9.97$E$-07$	&	44486	&	63899	&	834.13 	&	$1.00$E$-06$	\\
{\spalbb}(2)	&	799	&	41186	&	103.26 	&	$9.99$E$-07$	&	3052	&	37812	&	\pmb{490.00} 	&	$9.97$E$-07$	\\
{\spalbb}(3)	&	152	&	51307	&	133.36 	&	$6.43$E$-07$	&	429	&	72440	&	930.67 	&	$9.88$E$-07$	\\
{\spalbb}(4)	&	25	&	33030	&	86.80 	&	$1.00$E$-06$	&	56	&	51184	&	667.27 	&	$9.53$E$-07$	\\
{\spalbb}(5)	&	16	&	51574	&	139.75 	&	$1.00$E$-06$	&	16	&	54033	&	707.80 	&	$1.00$E$-06$	\\\hline
\end{tabular}
\end{table}

\begin{table}[htbp]
\setlength\tabcolsep{4.6pt}
\centering\scriptsize
\caption{Numerical results for \Cref{ex3} with $\nu=0.01$. }\label{tab8}
\begin{tabular}{|c|r|r|r|r|r|r|r|r|r|r|r|} \hline
$h (n,m)$&\multicolumn{4}{c|}{$2^{-5}$ $(n=5,890,m=769)$}
&\multicolumn{4}{c|}{$2^{-6}$ $(n=23,042,m=2,945)$}\\\cline{2-9}
  &Oiter&Titer&CPU&RES&Oiter&Titer&CPU&RES\\\hline
BICGSTAB	&		&	-	&	-	&	-	&		&	4746.5	&	\pmb{4.53} 	&	$9.85$E$-07$	\\
GMRES(20)	&		&	5518	&	2.17 	&	$1.00$E$-06$	&		&	10529	&	14.70 	&	$1.00$E$-06$	\\
GMRES(50)	&		&	3782	&	2.87 	&	$1.00$E$-06$	&		&	9379	&	26.42 	&	$1.00$E$-06$	\\
{\spalbb}(1)	&	584	&	8002	&	\pmb{0.73} 	&	$9.96$E$-07$	&	2340	&	9722	&	4.70 	&	$1.00$E$-06$	\\
{\spalbb}(2)	&	106	&	12810	&	1.15 	&	$9.99$E$-07$	&	297	&	20437	&	9.54 	&	$8.66$E$-07$	\\
{\spalbb}(3)	&	24	&	10581	&	0.91 	&	$4.68$E$-07$	&	53	&	17322	&	8.02 	&	$9.98$E$-07$	\\
{\spalbb}(4)	&	17	&	16101	&	1.42 	&	$9.99$E$-07$	&	18	&	18021	&	8.45 	&	$9.99$E$-07$	\\
{\spalbb}(5)	&	18	&	19627	&	1.75 	&	$1.00$E$-06$	&	17	&	26925	&	13.12 	&	$1.00$E$-06$	\\
\hline

$h (n,m)$&\multicolumn{4}{c|}{$2^{-7}$ $(n=91,138,m=11,521)$}
&\multicolumn{4}{c|}{$2^{-8}$ $(n=362,498,m=45,569)$}\\\cline{2-9}
  &Oiter&Titer&CPU&RES&Oiter&Titer&CPU&RES\\\hline
BICGSTAB	&		&	9904.5	&	48.44 	&	$9.99$E$-07$	&		&	-	&	-	&	-	\\
GMRES(20)	&		&	19010	&	104.40 	&	$9.99$E$-07$	&		&	37453	&	841.47 	&	$1.00$E$-06$	\\
GMRES(50)	&		&	19738	&	187.88 	&	$1.00$E$-06$	&		&	37074	&	1319.69 	&	$1.00$E$-06$	\\
{\spalbb}(1)	&	9286	&	19210	&	\pmb{45.65} 	&	$9.98$E$-07$	&	34222	&	47839	&	825.49 	&	$9.97$E$-07$	\\
{\spalbb}(2)	&	781	&	24365	&	61.89 	&	$1.00$E$-06$	&	3013	&	28458	&	\pmb{343.04} 	&	$1.00$E$-06$	\\
{\spalbb}(3)	&	145	&	29784	&	77.31 	&	$9.74$E$-07$	&	375	&	46354	&	557.08 	&	$9.99$E$-07$	\\
{\spalbb}(4)	&	25	&	23353	&	62.19 	&	$9.23$E$-07$	&	60	&	41959	&	500.22 	&	$9.99$E$-07$	\\
{\spalbb}(5)	&	16	&	39604	&	103.42 	&	$1.00$E$-06$	&	17	&	41434	&	507.15 	&	$1.00$E$-06$	\\
\hline
\end{tabular}
\end{table}

\begin{table}[htbp]
\setlength\tabcolsep{4.6pt}
\centering\scriptsize
\caption{Numerical results for \Cref{ex3} with $\nu=0.05$. }\label{tab9}
\begin{tabular}{|c|r|r|r|r|r|r|r|r|r|r|r|} \hline
$h (n,m)$&\multicolumn{4}{c|}{$2^{-5}$ $(n=5,890,m=769)$}
&\multicolumn{4}{c|}{$2^{-6}$ $(n=23,042,m=2,945)$}\\\cline{2-9}
  &Oiter&Titer&CPU&RES&Oiter&Titer&CPU&RES\\\hline
BICGSTAB	&		&	914.5	&	\pmb{0.22} 	&	$8.68$E$-07$	&		&	1888.5	&	\pmb{1.89} 	&	$2.44$E$-07$	\\
GMRES(20)	&		&	3139	&	1.35 	&	$9.98$E$-07$	&		&	6106	&	8.44 	&	$9.99$E$-07$	\\
GMRES(50)	&		&	2808	&	2.14 	&	$9.99$E$-07$	&		&	6467	&	16.24 	&	$9.99$E$-07$	\\
{\spalbb}(1)	&	427	&	2487	&	0.27 	&	$9.81$E$-07$	&	1449	&	5128	&	2.32 	&	$1.00$E$-06$	\\
{\spalbb}(2)	&	81	&	3657	&	0.33 	&	$8.61$E$-07$	&	196	&	6386	&	3.11 	&	$9.99$E$-07$	\\
{\spalbb}(3)	&	24	&	4135	&	0.37 	&	$1.00$E$-06$	&	45	&	7116	&	3.23 	&	$7.23$E$-07$	\\
{\spalbb}(4)	&	19	&	6332	&	0.59 	&	$1.00$E$-06$	&	19	&	8819	&	4.27 	&	$9.98$E$-07$	\\
{\spalbb}(5)	&	18	&	7035	&	0.65 	&	$9.98$E$-07$	&	18	&	11447	&	5.33 	&	$1.00$E$-06$	\\
\hline

$h (n,m)$&\multicolumn{4}{c|}{$2^{-7}$ $(n=91,138,m=11,521)$}
&\multicolumn{4}{c|}{$2^{-8}$ $(n=362,498,m=45,569)$}\\\cline{2-9}
  &Oiter&Titer&CPU&RES&Oiter&Titer&CPU&RES\\\hline
BICGSTAB	&		&	3811.5	&	\pmb{18.29} 	&	$6.50$E$-07$	&		&	8153.5	&	\pmb{194.33} 	&	$8.93$E$-07$	\\
GMRES(20)	&		&	12413	&	64.26 	&	$9.99$E$-07$	&		&	23112	&	526.24 	&	$1.00$E$-06$	\\
GMRES(50)	&		&	11779	&	99.11 	&	$9.99$E$-07$	&		&	21620	&	1401.43 	&	$1.00$E$-06$	\\
{\spalbb}(1)	&	4898	&	13457	&	30.89 	&	$9.98$E$-07$	&	15210	&	36389	&	441.11 	&	$9.99$E$-07$	\\
{\spalbb}(2)	&	542	&	9821	&	22.50 	&	$9.97$E$-07$	&	1601	&	21727	&	264.76 	&	$1.00$E$-06$	\\
{\spalbb}(3)	&	101	&	12949	&	30.46 	&	$9.98$E$-07$	&	250	&	21972	&	267.31 	&	$1.00$E$-06$	\\
{\spalbb}(4)	&	27	&	13639	&	33.74 	&	$6.40$E$-07$	&	46	&	22431	&	274.17 	&	$9.98$E$-07$	\\
{\spalbb}(5)	&	18	&	17621	&	44.04 	&	$9.99$E$-07$	&	19	&	26577	&	328.06 	&	$1.00$E$-06$	\\
\hline
\end{tabular}
\end{table}


\begin{ex}\label{ex:sd}
    Fluid flow in $\Omega_{f}\subset \R^2$ coupled with porous media flow in $\Omega_{p}\subset \R^2$ is governed by the static Stokes equations
\begin{equation}\label{sdm}
      -\nu\Delta\,{\bm u}_f + \nabla\, p_f  = {\bm f},
      \quad \textup{and} \quad
      {\rm div}\,{\bm u}_f = 0,
      \quad {\bm z}\in \Omega_{f},
\end{equation}
where $\Omega_{f}\cap \Omega_{p}=\varnothing$ and $\overline{\Omega}_{f}\cap \overline{\Omega}_{p}=\Gamma$ with $\Gamma$ being an interface, $\nu>0$ is the kinematic viscosity, and $\bm f$ is the external force. In the porous media region, the governing variable is $\phi = \frac{p_p}{\rho_f g}$, where $p_p$ is the pressure in $\Omega_{p}$, $\rho_f$ is the fluid density, and $g$ is the acceleration due to gravity. The velocity ${\bm u}_p$ of the porous media flow is related to $\phi$ via Darcy's law and is also divergence free:
\begin{equation}\label{sdm2}
    {\bm u}_p = -\dfrac{\epsilon^2}{r\nu}\nabla \phi
      \quad \textup{and} \quad
      -{\rm div}\,{\bm u}_p = 0,
      \quad {\bm z}\in \Omega_{p},
\end{equation}
where $r$ is the volumetric porosity and $\epsilon$ the characteristic length of the porous media.
\end{ex}

In our numerical experiments, the computational domain is $\Omega_f=(0,1)\times(1,2)$, $\Omega_p=(0,1)\times(0,1)$ and the interface is $\Gamma=(0,1)\times\{1\}$. We use a uniform mesh with grid parameters $h=2^{-5},\,2^{-6},\,2^{-7},\,2^{-8}$ to decompose $\Omega_f$, P2--P1 elements in the fluid region, and P2 Lagrange elements in the porous media region. We set $r=1$ and $\epsilon=\sqrt{0.1\nu}$, and again test $\nu=0.005$, $0.01$, $0.05$. Applying finite element discretization to the mixed Stokes-Darcy model~\eqref{sdm}--\eqref{sdm2} with the Dirichlet no-flow boundary conditions leads to linear systems of form~\eqref{saddle1} with 
   $G=\pmat{ G_{11} & G_{12}\\
                    -G_{12}^T & \nu G_{22}} $
\citep{Cai2009}.
Here $G$ is UPD and $B$ has full column rank. The numerical results are reported in \Cref{tab10,tab11,tab12} and \Cref{fig:ns6-sd1}. According to \Cref{tab10,tab11,tab12}, all methods again perform better for larger $\nu$, and BICGSTAB requires the least CPU time in most cases, while {\spalbb} is more competitive for smaller $\nu$. For \Cref{ex:sd}, {\spalbb} prefers smaller $\omega$, such as $\omega=10^{-5}$.

\begin{table}[htbp]
\setlength\tabcolsep{4.0pt}
\centering\scriptsize
\caption{Numerical results for \Cref{ex:sd} with $\nu=0.005$. }\label{tab10}
\begin{tabular}{|c|r|r|r|r|r|r|r|r|r|r|r|} \hline
$h (n,m)$&\multicolumn{4}{c|}{$2^{-5}$ $(n=12,675,m=1,089)$}
&\multicolumn{4}{c|}{$2^{-6}$ $(n=49,923,m=4,225)$}\\\cline{2-9}
  &Oiter&Titer&CPU&RES&Oiter&Titer&CPU&RES\\\hline
BICGSTAB	&		&	3286.5	&	\pmb{0.96} 	&	$9.25$E$-07$	&		&	7366.5	&	9.22 	&	$9.90$E$-07$	\\				
GMRES(20)	&		&	7668	&	5.25 	&	$1.00$E$-06$	&		&	21422	&	50.66 	&	$1.00$E$-06$	\\				
GMRES(50)	&		&	6373	&	8.74 	&	$1.00$E$-06$	&		&	15042	&	66.02 	&	$1.00$E$-06$	\\				
{\spalbb}(1)	&	2048	&	9709	&	2.51 	&	$9.99$E$-07$	&	6972	&	29955	&	56.13 	&	$1.00$E$-06$	\\				
{\spalbb}(2)	&	230	&	7387	&	1.13 	&	$9.98$E$-07$	&	740	&	15065	&	13.15 	&	$9.98$E$-07$	\\				
{\spalbb}(3)	&	38	&	7709	&	1.08 	&	$9.98$E$-07$	&	91	&	14640	&	9.56 	&	$1.00$E$-06$	\\				
{\spalbb}(4)	&	18	&	8204	&	1.13 	&	$9.99$E$-07$	&	23	&	13798	&	\pmb{8.60} 	&	$1.00$E$-06$	\\				
{\spalbb}(5)	&	18	&	11310	&	1.77 	&	$9.99$E$-07$	&	18	&	25189	&	15.77 	&	$1.00$E$-06$	\\				
\hline
	
$h (n,m)$&\multicolumn{4}{c|}{$2^{-7}$ $(n=198,147,m=16,641)$}
&\multicolumn{4}{c|}{$2^{-8}$ $(n=789,507,m=66,049)$}\\\cline{2-9}
  &Oiter&Titer&CPU&RES&Oiter&Titer&CPU&RES\\\hline
BICGSTAB	&		&	15174.5	&	\pmb{107.22} 	&	$8.93$E$-07$	&		&	32027.5	&	829.92 	&	$9.79$E$-07$	\\				
GMRES(20)	&		&	47298	&	450.74 	&	$1.00$E$-06$	&		&	109310	&	4099.15 	&	$1.00$E$-06$	\\				
GMRES(50)	&		&	40626	&	605.13 	&	$1.00$E$-06$	&		&	81406	&	5300.94 	&	$1.00$E$-06$	\\				
{\spalbb}(1)	&	23244	&	94747	&	1623.05 	&	$1.00$E$-06$	&	-	&	-	&	-	&	-	\\				
{\spalbb}(2)	&	2442	&	43102	&	305.56 	&	$1.00$E$-06$	&	-	&	-	&	-	&	-	\\				
{\spalbb}(3)	&	272	&	30164	&	132.12 	&	$1.00$E$-06$	&	815	&	66624	&	1503.82 	&	$9.99$E$-07$	\\				
{\spalbb}(4)	&	42	&	27889	&	113.17 	&	$9.99$E$-07$	&	98	&	56613	&	907.24 	&	$9.99$E$-07$	\\				
{\spalbb}(5)	&	18	&	33198	&	138.39 	&	$1.00$E$-06$	&	24	&	47077	&	\pmb{717.58} 	&	$1.00$E$-06$	\\				
\hline
\end{tabular}
\end{table}

\begin{table}[htbp]
\setlength\tabcolsep{4.2pt}
\centering\scriptsize
\caption{Numerical results for \Cref{ex:sd}  with $\nu=0.01$. }\label{tab11}
\begin{tabular}{|c|r|r|r|r|r|r|r|r|r|r|r|} \hline
$h (n,m)$&\multicolumn{4}{c|}{$2^{-5}$ $(n=12,675,m=1,089)$}
&\multicolumn{4}{c|}{$2^{-6}$ $(n=49,923,m=4,225)$}\\\cline{2-9}
  &Oiter&Titer&CPU&RES&Oiter&Titer&CPU&RES\\\hline
BICGSTAB	& &	2485.5	&	\pmb{0.66} 	&	$9.88$E$-07$&&	4888.5	&	\pmb{5.27} 	&	$9.95$E$-07$	\\
GMRES(20)	& &	5717	&	3.44 	&	$9.99$E$-07$&&	14398	&	29.13 	&	$9.99$E$-07$	 \\
GMRES(50)	& &	4276	&	5.27 	&	$1.00$E$-06$&&	11608	&	47.31 	&	$1.00$E$-06$		\\
{\spalbb}(1)	&	2346	&	9426	&	2.59 	&	$9.99$E$-07$&	7762	&	30273	&	58.77 	&	$1.00$E$-06$	\\
{\spalbb}(2)	&	265	&	5132	&	0.78 	&	$9.99$E$-07$&	829	&	13805	&	12.72 	&	$1.00$E$-06$	\\
{\spalbb}(3)	&	43	&	5380	&	0.69 	&	$9.97$E$-07$&	104	&	11902	&	7.44 	&	$9.99$E$-07$	\\
{\spalbb}(4)	&	18	&	6560	&	0.85 	&	$9.98$E$-07$&	25	&	11578	&	7.05 	&	$9.99$E$-07$	\\
{\spalbb}(5)	&	18	&	9092	&	1.14 	&	$1.00$E$-06$&	18	&	17747	&	10.45 	&	$9.98$E$-07$	\\\hline

$h (n,m)$&\multicolumn{4}{c|}{$2^{-7}$ $(n=198,147,m=16,641)$}
&\multicolumn{4}{c|}{$2^{-8}$ $(n=789,507,m=66,049)$}\\\cline{2-9}
  &Oiter&Titer&CPU&RES&Oiter&Titer&CPU&RES\\\hline
BICGSTAB	& &	10216	&	\pmb{70.41} 	&	$9.78$E$-07$&&	20513.5	&	\pmb{579.78} 	&	$9.94$E$-07$	\\
GMRES(20)	&		&	33336	&	326.81 	&	$1.00$E$-06$&&	58210	&	6339.92 	&	$1.00$E$-06$	\\
GMRES(50)	& &	28362	&	438.14 	&	$1.00$E$-06$&&	48116	&	3648.93 	&	$1.00$E$-06$	\\
{\spalbb}(1)	&	25411	&	97928	&	1769.82 	&	$1.00$E$-06$&-	&	-	&	- 	&	- \\
{\spalbb}(2)	&	2680	&	42275	&	318.61 	&	$9.99$E$-07$&-	&	-	&	- 	&	- \\
{\spalbb}(3)	&	295	&	23860	&	111.65 	&	$9.97$E$-07$&877	&	61065	&	1477.52 	&	$9.98$E$-07$	\\
{\spalbb}(4)	&	46	&	22297	&	93.48 	&	$9.98$E$-07$&109	&	43631	&	703.09 	&	$1.00$E$-06$	\\
{\spalbb}(5)	&	19	&	21946	&	92.28 	&	$1.00$E$-06$&26	&	37877	&	583.34 	&	$9.99$E$-07$	\\		
\hline
\end{tabular}
\end{table}

\begin{table}[htbp]
\setlength\tabcolsep{4.6pt}
\centering\scriptsize
\caption{Numerical results for \Cref{ex:sd}  with $\nu=0.05$. }\label{tab12}
\begin{tabular}{|c|r|r|r|r|r|r|r|r|r|r|r|} \hline
$h (n,m)$&\multicolumn{4}{c|}{$2^{-5}$ $(n=12,675,m=1,089)$}
&\multicolumn{4}{c|}{$2^{-6}$ $(n=49,923,m=4,225)$}\\\cline{2-9}
  &Oiter&Titer&CPU&RES&Oiter&Titer&CPU&RES\\\hline
BICGSTAB	&		&	883.5	&	\pmb{0.22} 	&	$9.94$E$-07$	&		&	1738.5	&	\pmb{1.57} 	&	$9.55$E$-07$	\\
GMRES(20)	&		&	2589	&	2.22 	&	$9.98$E$-07$	&		&	5115	&	9.27 	&	$9.99$E$-07$	\\
GMRES(50)	&		&	2509	&	3.85 	&	$1.00$E$-06$	&		&	4542	&	16.25 	&	$1.00$E$-06$	\\
{\spalbb}(1)	&	3077	&	10839	&	2.94 	&	$9.97$E$-07$	&	10033	&	35063	&	68.21 	&	$1.00$E$-06$	\\
{\spalbb}(2)	&	346	&	4649	&	0.70 	&	$9.99$E$-07$	&	1090	&	14684	&	12.95 	&	$9.98$E$-07$	\\
{\spalbb}(3)	&	54	&	3899	&	0.46 	&	$9.98$E$-07$	&	131	&	7892	&	4.30 	&	$9.94$E$-07$	\\
{\spalbb}(4)	&	21	&	3344	&	0.38 	&	$9.99$E$-07$	&	29	&	6819	&	3.38 	&	$9.98$E$-07$	\\
{\spalbb}(5)	&	18	&	3263	&	0.37 	&	$9.97$E$-07$	&	18	&	6022	&	2.95 	&	$1.00$E$-06$	\\

\hline

$h (n,m)$&\multicolumn{4}{c|}{$2^{-7}$ $(n=198,147,m=16,641)$}
&\multicolumn{4}{c|}{$2^{-8}$ $(n=789,507,m=66,049)$}\\\cline{2-9}
  &Oiter&Titer&CPU&RES&Oiter&Titer&CPU&RES\\\hline
BICGSTAB	&		&	3512.5	&	\pmb{23.21} 	&	$9.92$E$-07$	&&	7305.5	&	\pmb{230.70} 	&	$9.72$E$-07$	\\
GMRES(20)	&		&	10545	&	107.69 	&	$1.00$E$-06$	&&	30657	&	1438.72 	&	$1.00$E$-06$	\\
GMRES(50)	&		&	8036	&	164.70 	&	$1.00$E$-06$	&&	14974	&	1230.40 	&	$1.00$E$-06$	\\
{\spalbb}(1)	&	-	&	-	&	- 	&	-	& -	&	-	&	- 	&	-	\\
{\spalbb}(2)	&	3464	&	46071	&	375.30 	&	$9.98$E$-07$	&-	&	-	&	-	&	-	\\
{\spalbb}(3)	&	373	&	21106	&	102.87 	&	$9.97$E$-07$	&1105	&	60707	&	1291.67 	&	$9.96$E$-07$\\
{\spalbb}(4)	&	55	&	13647	&	56.56 	&	$9.99$E$-07$	&	126	&	26414	&	469.18 	&	$1.00$E$-06$\\
{\spalbb}(5)	&	20	&	11646	&	48.65 	&	$1.00$E$-06$	&   28	&	21844	&	341.45 	&	$9.99$E$-07$\\

\hline
\end{tabular}
\end{table}


In conclusion, \Cref{tab1,tab2,tab3,tab4,tab5,tab6,tab7,tab8,tab9,tab10,tab11,tab12} and \Cref{fig:ns3-ns5,fig:ns6-sd1} illustrate that {\spalbb} is a practical method, and its advantages increase with problem size. {\spalbb} and GMRES are more robust than BICGSTAB. Unlike GMRES, {\spalbb} has constant storage. In terms of CPU time, {\spalbb} is more efficient than GMRES. We see from \Cref{tab1,tab2,tab3,tab4,tab5,tab6,tab7,tab8,tab9} that the advantages of {\spalbb} are more obvious for smaller $\nu$, i.e., more unsymmetric $G$. \Cref{fig:ns3-ns5,fig:ns6-sd1} indicate that the convergence rate of {\spalbb} depends strongly on $\omega$. For larger $\omega$, the nonmonotonicity of $\|r_k\|$ in {\spalbb} becomes more pronounced. The strong nonmonotone behavior is similar to the BB method \cite{raydan1997}.

\section{Conclusions}\label{sec:con}
We presented a theoretical and numerical study of the augmented Lagrangian ({\spal}) algorithm and its inexact version for solving unsymmetric saddle-point systems. Specifically, we used a gradient method, known as the Barzilai-Borwein (BB) method, to solve the linear system in {\spal} inexactly and proposed the augmented Lagrangian BB ({\spalbb}) algorithm. The numerical results for {\spalbb} presented are highly encouraging. {\spalbb} often requires the least CPU time, and, especially for larger problems, its advantages are clear. Practical methods for choosing $\omega$ and $Q$ to balance the inner and outer iterations is a topic for future research. 

\section*{Acknowledgments}
We thank our colleague and friend, Prof Dr Oleg Burdakov, for his devotion to this research. In particular, we express our gratitude to him for fundamental contributions that initiated this work. Oleg developed the {\spalbb} algorithm, proposed the counter-example to show that the BB1 method may be divergent, and gave many constructive suggestions on our Matlab implementation of {\spalbb}.

\small
\bibliographystyle{abbrvnat}
\bibliography{references}

\begin{thebibliography}{52}
\providecommand{\natexlab}[1]{#1}
\providecommand{\url}[1]{\texttt{#1}}
\expandafter\ifx\csname urlstyle\endcsname\relax
  \providecommand{\doi}[1]{doi: #1}\else
  \providecommand{\doi}{doi: \begingroup \urlstyle{rm}\Url}\fi

\bibitem[Arrow et~al.(1958)Arrow, Hurwicz, Uzawa, and Chenery]{arrow1958}
K.~J. Arrow, L.~Hurwicz, H.~Uzawa, and H.~B. Chenery.
\newblock \emph{{Studies in Linear and Non-linear Programming}}, volume~2.
\newblock Stanford University Press, 1958.

\bibitem[Awanou and Lai(2005{\natexlab{a}})]{awanou2005}
G.~Awanou and M.~J. Lai.
\newblock On convergence rate of the augmented {Lagrangian} algorithm for
  nonsymmetric saddle point problems.
\newblock \emph{Appl. Numer. Math.}, 54\penalty0 (2):\penalty0 122--134,
  2005{\natexlab{a}}.

\bibitem[Awanou and Lai(2005{\natexlab{b}})]{awanou2005trivariate}
G.~Awanou and M.~J. Lai.
\newblock {Trivariate spline approximations of 3D Navier--Stokes equations}.
\newblock \emph{Math. Comp.}, 74\penalty0 (250):\penalty0 585--601,
  2005{\natexlab{b}}.

\bibitem[Bai and Benzi(2017)]{bai2017regularized}
Z.~Z. Bai and M.~Benzi.
\newblock Regularized {HSS} iteration methods for saddle-point linear systems.
\newblock \emph{BIT Numer. Math.}, 57\penalty0 (2):\penalty0 287--311, 2017.

\bibitem[Barzilai and Borwein(1988)]{Barzilai1988two}
J.~Barzilai and J.~M. Borwein.
\newblock Two-point step size gradient methods.
\newblock \emph{IMA J. Numer. Anal.}, 8\penalty0 (1):\penalty0 141--148, 1988.

\bibitem[Benzi and Golub(2004)]{benzi2004}
M.~Benzi and G.~H. Golub.
\newblock A preconditioner for generalized saddle point problems.
\newblock \emph{SIAM J. Matrix Anal. Appl.}, 26\penalty0 (1):\penalty0 20--41,
  2004.

\bibitem[Benzi and Wathen(2008)]{benzi2008some}
M.~Benzi and A.~J. Wathen.
\newblock Some preconditioning techniques for saddle point problems.
\newblock \emph{Model order reduction: theory, research aspects and
  applications}, pages 195--211, 2008.

\bibitem[Benzi et~al.(2005)Benzi, Golub, and Liesen]{Benzi2005}
M.~Benzi, G.~H. Golub, and J.~Liesen.
\newblock Numerical solution of saddle point problems.
\newblock \emph{Acta Numerica}, 14\penalty0 (2):\penalty0 1--137, 2005.

\bibitem[Berman and Plemmons(1994)]{Berman1994}
A.~Berman and R.~J. Plemmons.
\newblock \emph{{Nonnegative Matrices in the Mathematical Sciences}}.
\newblock SIAM, 1994.

\bibitem[Bertsekas(2014)]{bertsekas2014}
D.~P. Bertsekas.
\newblock \emph{{Constrained Optimization and Lagrange Multiplier Methods}}.
\newblock Academic Press, 2014.

\bibitem[Birgin and Mart{\'\i}nez(2014)]{birgin2014}
E.~G. Birgin and J.~M. Mart{\'\i}nez.
\newblock \emph{{Practical Augmented Lagrangian Methods for Constrained
  Optimization}}.
\newblock SIAM, 2014.

\bibitem[Burdakov et~al.(2019)Burdakov, Dai, and Huang]{burdakov2019}
O.~Burdakov, Y.~H. Dai, and N.~Huang.
\newblock {Stabilized Barzilai-Borwein method}.
\newblock \emph{J. Comput. Math.}, 37\penalty0 (6):\penalty0 916--936, 2019.

\bibitem[Cai et~al.(2009)Cai, Mu, and Xu]{Cai2009}
M.~C. Cai, M.~Mu, and J.~C. Xu.
\newblock Preconditioning techniques for a mixed {Stokes/Darcy} model in porous
  media applications.
\newblock \emph{J. Comput. Appl. Math.}, 233\penalty0 (2):\penalty0 346--355,
  2009.

\bibitem[Campbell and Meyer(2009)]{Campbell1979}
S.~L. Campbell and C.~D. Meyer.
\newblock \emph{{Generalized Inverses of Linear Transformations}}.
\newblock SIAM, 2009.

\bibitem[Cao and Miao(2016)]{cao2016}
Y.~Cao and S.~X. Miao.
\newblock On semi-convergence of the generalized shift-splitting iteration
  method for singular nonsymmetric saddle point problems.
\newblock \emph{Comput. Math. Appl.}, 71\penalty0 (7):\penalty0 1503--1511,
  2016.

\bibitem[Cheng(2000)]{cheng2000}
X.~L. Cheng.
\newblock On the nonlinear inexact {Uzawa} algorithm for saddle-point problems.
\newblock \emph{SIAM J. Numer. Anal.}, 37\penalty0 (6):\penalty0 1930--1934,
  2000.

\bibitem[Dai and Liao(2002)]{dai2002r}
Y.~H. Dai and L.~Z. Liao.
\newblock {R-linear convergence of the Barzilai and Borwein gradient method}.
\newblock \emph{IMA J. Numer. Anal.}, 22\penalty0 (1):\penalty0 1--10, 2002.

\bibitem[Dai et~al.(2005)Dai, Liao, and Li]{dai2005analysis}
Y.~H. Dai, L.~Z. Liao, and D.~Li.
\newblock {An analysis of Barzilai-Borwein gradient method for unsymmetric
  linear equations}.
\newblock \emph{Optim. Control Appl.}, pages 183--211, 2005.

\bibitem[Dai et~al.(2006)Dai, Hager, Schittkowski, and Zhang]{dai2006cyclic}
Y.~H. Dai, W.~W. Hager, K.~Schittkowski, and H.~C. Zhang.
\newblock {The cyclic Barzilai-Borwein method for unconstrained optimization}.
\newblock \emph{IMA J. Numer. Anal.}, 26\penalty0 (3):\penalty0 604--627, 2006.

\bibitem[Di~Serafino and Orban(2021)]{di2021constraint}
D.~Di~Serafino and D.~Orban.
\newblock Constraint-preconditioned {Krylov} solvers for regularized
  saddle-point systems.
\newblock \emph{SIAM J. Sci. Comput.}, 43\penalty0 (2):\penalty0 A1001--A1026,
  2021.

\bibitem[Dollar et~al.(2010)Dollar, Gould, Stoll, and Wathen]{dollar2010}
H.~S. Dollar, N.~I. Gould, M.~Stoll, and A.~J. Wathen.
\newblock Preconditioning saddle-point systems with applications in
  optimization.
\newblock \emph{SIAM J. Sci. Comput.}, 32\penalty0 (1):\penalty0 249--270,
  2010.

\bibitem[Elman(1999)]{elman1999preconditioning}
H.~C. Elman.
\newblock Preconditioning for the steady-state {Navier--Stokes} equations with
  low viscosity.
\newblock \emph{SIAM Journal on Scientific Computing}, 20\penalty0
  (4):\penalty0 1299--1316, 1999.

\bibitem[Elman et~al.(2007)Elman, Ramage, and Silvester]{elman2007}
H.~C. Elman, A.~Ramage, and D.~J. Silvester.
\newblock {Algorithm 866: IFISS, a Matlab toolbox for modelling incompressible
  flow}.
\newblock \emph{ACM Trans. Math. Softw.}, 33\penalty0 (2):\penalty0 14--es,
  2007.

\bibitem[Friedlander et~al.(1998)Friedlander, Mart\'{\i}nez, Molina, and
  Raydan]{friedlander1998}
A.~Friedlander, J.~M. Mart\'{\i}nez, B.~Molina, and M.~Raydan.
\newblock Gradient method with retards and generalizations.
\newblock \emph{SIAM J. Numer. Anal.}, 36\penalty0 (1):\penalty0 275--289,
  1998.

\bibitem[Ghannad et~al.(2022)Ghannad, Orban, and
  Saunders]{ghannad-orban-saunders-2021}
A.~Ghannad, D.~Orban, and M.~A. Saunders.
\newblock \doilink{10.1080/10556788.2021.1965599}{Linear systems arising in
  interior methods for convex optimization: a symmetric formulation with
  bounded condition number}.
\newblock \emph{Optim. Method Softw.}, 37\penalty0 (4):\penalty0 1344--1369,
  2022.

\bibitem[Glowinski and Le~Tallec(1989)]{glowinski1989}
R.~Glowinski and P.~Le~Tallec.
\newblock \emph{{Augmented Lagrangian and Operator-splitting Methods in
  Nonlinear Mechanics}}.
\newblock SIAM, 1989.

\bibitem[Golub and Greif(2003)]{golub2003solving}
G.~H. Golub and C.~Greif.
\newblock On solving block-structured indefinite linear systems.
\newblock \emph{SIAM J. Sci. Comput.}, 24\penalty0 (6):\penalty0 2076--2092,
  2003.

\bibitem[Golub et~al.(2005)Golub, Greif, and Varah]{greif2004}
G.~H. Golub, C.~Greif, and J.~M. Varah.
\newblock An algebraic analysis of a block diagonal preconditioner for saddle
  point systems.
\newblock \emph{SIAM J. Matrix Anal. Appl.}, 27\penalty0 (3):\penalty0
  779--792, 2005.

\bibitem[Gould et~al.(2014)Gould, Orban, and Rees]{gould2014}
N.~Gould, D.~Orban, and T.~Rees.
\newblock Projected {Krylov} methods for saddle-point systems.
\newblock \emph{SIAM J. Matrix Anal. Appl.}, 35\penalty0 (4):\penalty0
  1329--1343, 2014.

\bibitem[Hu and Zou(2006)]{hu2006nonlinear}
Q.~Hu and J.~Zou.
\newblock Nonlinear inexact {Uzawa} algorithms for linear and nonlinear
  saddle-point problems.
\newblock \emph{SIAM J. Optim.}, 16\penalty0 (3):\penalty0 798--825, 2006.

\bibitem[Kozjakin and Krasnosel'ski(1982)]{kozjakin1982some}
V.~Kozjakin and M.~Krasnosel'ski.
\newblock Some remarks on the method of minimal residues.
\newblock \emph{Numer. Funct. Anal. Optim.}, 4\penalty0 (3):\penalty0 211--239,
  1982.

\bibitem[Krasnosel'skii and Krein(1952)]{krasnosel1952}
M.~A. Krasnosel'skii and S.~G. Krein.
\newblock An iteration process with minimal residuals.
\newblock \emph{Matematicheskii Sbornik}, 73\penalty0 (2):\penalty0 315--334,
  1952.

\bibitem[Lu and Zhang(2010)]{lu2010}
J.~Lu and Z.~Zhang.
\newblock A modified nonlinear inexact {Uzawa} algorithm with a variable
  relaxation parameter for the stabilized saddle point problem.
\newblock \emph{SIAM J. Matrix Anal. Appl.}, 31\penalty0 (4):\penalty0
  1934--1957, 2010.

\bibitem[Molina and Raydan(1996)]{molina1996}
B.~Molina and M.~Raydan.
\newblock Preconditioned {Barzilai-Borwein} method for the numerical solution
  of partial differential equations.
\newblock \emph{Numer. Algor.}, 13:\penalty0 45--60, 1996.

\bibitem[Montoison and Orban(2023)]{montoison2023gpmr}
A.~Montoison and D.~Orban.
\newblock {GPMR}: An iterative method for unsymmetric partitioned linear
  systems.
\newblock \emph{SIAM J. Matrix Anal. Appl.}, 44\penalty0 (1):\penalty0
  293--311, 2023.

\bibitem[Orban and Arioli(2017)]{orban2017}
D.~Orban and M.~Arioli.
\newblock {Full-Space Iterative Methods}.
\newblock In \emph{{Iterative Solution of Symmetric Quasi-definite Linear
  Systems}}, chapter~6, pages 63--72. SIAM, 2017.

\bibitem[Pestana and Rees(2016)]{pestana2016null}
J.~Pestana and T.~Rees.
\newblock Null-space preconditioners for saddle point systems.
\newblock \emph{SIAM J. Matrix Anal. Appl.}, 37\penalty0 (3):\penalty0
  1103--1128, 2016.

\bibitem[Ramage and Gartland~Jr(2013)]{Ramage2013}
A.~Ramage and E.~C. Gartland~Jr.
\newblock A preconditioned nullspace method for liquid crystal director
  modeling.
\newblock \emph{SIAM J. Sci. Comput.}, 35\penalty0 (1):\penalty0 B226--B247,
  2013.

\bibitem[Raydan(1993)]{raydan1993}
M.~Raydan.
\newblock {On the Barzilai and Borwein choice of steplength for the gradient
  method}.
\newblock \emph{IMA J. Numer. Anal.}, 13\penalty0 (3):\penalty0 321--326, 1993.

\bibitem[Raydan(1997)]{raydan1997}
M.~Raydan.
\newblock {The Barzilai and Borwein gradient method for the large scale
  unconstrained minimization problem}.
\newblock \emph{SIAM J. Optim.}, 7\penalty0 (1):\penalty0 26--33, 1997.

\bibitem[Rozlozn{\'\i}k and Simoncini(2002)]{rozloznik2002krylov}
M.~Rozlozn{\'\i}k and V.~Simoncini.
\newblock {Krylov} subspace methods for saddle point problems with indefinite
  preconditioning.
\newblock \emph{SIAM J. Matrix Anal. Appl.}, 24\penalty0 (2):\penalty0
  368--391, 2002.

\bibitem[Saad(2003)]{saad2003}
Y.~Saad.
\newblock \emph{{Iterative Methods for Sparse Linear Systems}}.
\newblock SIAM, 2003.

\bibitem[Saad and Schultz(1986)]{Saad1986}
Y.~Saad and M.~H. Schultz.
\newblock {GMRES}: A generalized minimal residual algorithm for solving
  nonsymmetric linear systems.
\newblock \emph{SIAM J. Sci. and Statist. Comput.}, 7\penalty0 (3):\penalty0
  856--869, 1986.

\bibitem[Scott and Tuma(2020)]{scott2020null}
J.~Scott and M.~Tuma.
\newblock A null-space approach for symmetric saddle point systems with a non
  zero (2, 2) block.
\newblock \emph{SIAM J. Sci. Comput.}, 2020.

\bibitem[Scott and Tuma(2022)]{scott2022null}
J.~Scott and M.~Tuma.
\newblock A null-space approach for large-scale symmetric saddle point systems
  with a small and non zero (2, 2) block.
\newblock \emph{Numer. Algor.}, 90\penalty0 (4):\penalty0 1639--1667, 2022.

\bibitem[Silvester et~al.(2023)Silvester, Elman, and Ramage]{ifiss}
D.~Silvester, H.~Elman, and A.~Ramage.
\newblock {Incompressible Flow \& Iterative Solver Software}.
\newblock
  \url{https://personalpages.manchester.ac.uk/staff/david.silvester/ifiss/},
  2023.

\bibitem[Van~der Vorst(1992)]{Van1992}
H.~A. Van~der Vorst.
\newblock {Bi-CGSTAB}: A fast and smoothly converging variant of {Bi-CG} for
  the solution of nonsymmetric linear systems.
\newblock \emph{SIAM J. Sci. and Statist. Comput.}, 13\penalty0 (2):\penalty0
  631--644, 1992.

\bibitem[Wright(1997)]{Wright1997}
S.~J. Wright.
\newblock \emph{{Primal-dual Interior-point Methods}}.
\newblock SIAM, 1997.

\bibitem[Zhang and Wei(2010)]{Zhang2010}
N.~Zhang and Y.~M. Wei.
\newblock On the convergence of general stationary iterative methods for
  {range-Hermitian} singular linear systems.
\newblock \emph{Numer. Linear Algebra Appl.}, 17:\penalty0 139--154, 2010.

\bibitem[Zheng et~al.(2009)Zheng, Bai, and Yang]{zheng2009}
B.~Zheng, Z.~Z. Bai, and X.~Yang.
\newblock On semi-convergence of parameterized {Uzawa} methods for singular
  saddle point problems.
\newblock \emph{Linear Algebra Appl.}, 431\penalty0 (5-7):\penalty0 808--817,
  2009.

\bibitem[Zou and Magoul\`{e}s(2022)]{zou2022delayed}
Q.~M. Zou and F.~Magoul\`{e}s.
\newblock Delayed gradient methods for symmetric and positive definite linear
  systems.
\newblock \emph{SIAM Rev.}, 64\penalty0 (3):\penalty0 517--553, 2022.

\bibitem[Zulehner(2002)]{zulehner2002analysis}
W.~Zulehner.
\newblock Analysis of iterative methods for saddle point problems: a unified
  approach.
\newblock \emph{Math. Comp.}, 71\penalty0 (238):\penalty0 479--505, 2002.

\end{thebibliography}
\normalsize

\end{document}